\documentclass[3p,numbers]{elsarticle}

\usepackage{graphicx}
\usepackage[section]{placeins}

\usepackage{subcaption}
\usepackage{booktabs} %
\usepackage{array} %
\usepackage{amsmath,amssymb,amsfonts}
\usepackage{amsthm, mathtools}

\usepackage[utf8]{inputenc}    %
\usepackage[T1]{fontenc}       %
\usepackage{lmodern}           %

\usepackage[UKenglish]{babel}    %
\usepackage[autostyle, english=british]{csquotes}
\usepackage[babel=true]{microtype}         %

\usepackage[breaklinks,unicode]{hyperref}
\usepackage[noabbrev,nameinlink]{cleveref}

\usepackage{xspace}

\newtheorem{thm}{Theorem}[section]
\newtheorem{lem}[thm]{Lemma}

\newtheorem{defn}[thm]{Definition}
\newtheorem{rmk}[thm]{Remark}

\newtheorem{exmp}[thm]{Example}

\newtheorem{meth}{Method}

\newcommand{\domain}{\Omega}
\newcommand{\domainbnd}{\partial\domain{}}
\newcommand{\extdomain}{B}
\newcommand{\extdomainbnd}{\partial\extdomain{}}
\newcommand{\DD}{Diffuse Domain}
\newcommand{\blfa}[2]{a(#1, #2)}
\newcommand{\blfl}[1]{l(#1)}

\newcommand{\BC}{\zeta(u, \phi, g)}
\newcommand{\tBC}{\tilde{\zeta}(u, \phi, g)}

\newcommand{\DDMO}{\mbox{DDM1}\xspace}
\newcommand{\MixedZ}{\mbox{Mix0DDM}\xspace}

\newcommand{\MixedO}{\mbox{Mix1DDM}\xspace}
\newcommand{\NSDDM}{\mbox{NSDDM}\xspace}
\newcommand{\NDDM}{\mbox{NDDM}\xspace}

\newcommand{\ddfem}{\mbox{\textsc{ddfem}}\xspace}
\newcommand{\dunefem}{\mbox{\textsc{Dune-Fem}}\xspace}
\newcommand{\dunealu}{\mbox{\textsc{Dune-ALUGrid}}\xspace}

\newcommand{\gmsh}{\mbox{Gmsh}\xspace}

\journal{Journal of Computational Physics}

\begin{document}
\begin{frontmatter}
    \title{Diffuse Domain Methods with Dirichlet Boundary Conditions}

    \author{Luke Benfield\corref{cor1}}
    \ead{luke.benfield@warwick.ac.uk}
    \cortext[cor1]{Corresponding author}

    \author{Andreas Dedner}
    \ead{a.s.dedner@warwick.ac.uk}

    \affiliation{organization={Mathematics Institute},
        addressline={University of Warwick},
        city={Coventry},
        postcode={CV4 7AL},
        country={UK}
    }

    \begin{abstract}

        The solution of partial differential equations (PDEs) on complex domains often
        presents a significant computational challenge by
        requiring the generation of fitted meshes.
        The Diffuse Domain Method (DDM) is an alternative which
        reformulates the problem on a larger, simple domain where
        the complex geometry is represented by a smooth phase-field function.

        This paper introduces and analyses several new DDM methods for
        solving problems with Dirichlet boundary conditions.
        We derive two new methods from the mixed formulation of the governing equations.
        This approach transforms the essential Dirichlet conditions into natural boundary conditions.
        Additionally, we develop coercive formulations based on Nitsche's
        method. We provide proofs of coercivity for all new methods and also
        for some approximations found in the literature for a wide range of linear
        advection-diffusion problems.

        Numerical experiments demonstrate the improved accuracy of the new methods,
        and reveal the balance between $L^2$ and $H^1$ errors.
        The practical effectiveness of this approach is demonstrated through
        the simulation of the incompressible Navier-Stokes equations on two
        fluid dynamics benchmark problems.

    \end{abstract}

\end{frontmatter}

\section{Introduction}
\label{section:intro2}

Solving partial differential equations (PDEs) on domains with complex geometries is
a significant challenge in computational science and engineering.
Traditional numerical methods usually depend on generating a fitted mesh,
which can be computationally expensive.
Due to the boundary complexities,
it is necessary to use a very fine mesh
that is often smaller than the solution's length scale.
This challenge is particularly present in three dimensions,
and problems involving evolving boundaries.
To overcome these difficulties of mesh generation,
a variety of methods have been developed.
For example, one may use a fitted mesh
and modifying the discretisation near the boundary,
such as the Extended Finite Element Method \citep{Fries2010}.
On the other hand, there are methods that utilise non-fitted meshes,
such as: the Fictitious Domain Method \citep{Glowinski1996,Ramiere2007,Vos2008,Parussini2009},
and the Immersed Boundary Method \citep{LeVeque1994,Prenter2023,Griffith2020}.

Among these approaches,
the \DD{} Method (DDM) was introduced by \citet{Kockelkoren2003} as a particularly simple approach.
The core concept is to embed the original, complex domain $\domain{}$
into a larger, simple domain $\extdomain{}$.
Therefore, it avoids the need for complex mesh generation, or explicit boundary tracking.
The problem's geometry is represented by
a phase-field function, $\phi$,
which smoothly transitions
between zero and one across a thin interfacial layer that approximates the boundary $\domainbnd{}$.
This transforms the problem of complex geometry into a problem of
modifying the system's equations with terms depending on $\phi$.
The main advantage of this method is its compatibility with standard finite element solvers,
as it does not require unique additional terms.

It is important to acknowledge that
the \DD{} method represents a collection of different approximations
which follow the above concept.
Asymptotic analysis has been shown for several approximations for
Dirichlet, Neumann, and Robin boundary conditions,
\citep{Li2009, Yu2020, Lervaag2015, Nguyen2017}.
These methods have been developed to work on surface domains \citep{Raetz2006},
and coupled bulk-surface systems \citep{Teigen2009, Abels2015}.
Also, there are similar approaches such as
the Smoothed Boundary method \citep{Yu2009, Yu2012, Termuhlen2022}.

The method has been successfully applied to a wide range of problems,
from biology \citep{Aland2012, Nguyen2017, Nguyen2017a, Chen2014}
to materials \citep{Chadwick2018, Raetz2016, Bukac2023}.
Our experimental goal will be to apply this to
the Navier-Stokes equations \citep{Aland2010, Guo2021, Anderson1998, Termuhlen2022}.
Consequently, one of our main objectives will be on second-order PDEs with
advection terms.

Typically, \DD{} method formulations are derived by modifying the strong form
or the weak form of the PDE and then including additional terms to incorporate boundary conditions.
This paper explores the development of new DDM formulations derived from
combining these together with the mixed formulation of the equations.
Hence, the Dirichlet boundary conditions become natural boundary conditions within the formulation,
which can lead to improved accuracy of the solution gradient.
Furthermore,
we follow the approach of Nitsche's method,
another well-established technique for the weak enforcement of boundary conditions \citep{Freund1995}
and has previously been used with the \DD{} method \citep{Nguyen2017, Monte2022}.
Along with the derivation of these new methods,
we will show proof of coercivity, showing how some existing methods need to
be modified to guarantee this important property of the approximation.

The paper is organised as follows.
\Cref{section:dd_background}
introduces the key concepts of the \DD{} method.
\Cref{section:Nitsche}
derives the Nitsche's based \DD{} approaches.
In \cref{sec:mixed_methods},
we introduce the new methods based on the mixed formulation.
\Cref{section:advection}
addresses the specific challenges and stabilisation terms
for advection-dominated problems.
\Cref{sec:Numerical}
presents a set of numerical results
focusing on application to advection-diffusion equations with Dirichlet boundary conditions,
and concluding with benchmark problems for the incompressible Navier-Stokes equations
to validate and compare the methods.

\section{Diffuse Domain Method}
\label{section:dd_background}

\subsection{Background}
\label{subsection:sharp2dd}

Starting with the example of Poisson's equation in a sufficiently smooth complex domain, $\domain{}$,
with Dirichlet boundary conditions:
\begin{subequations}
    \label{eq:poisson}
    \begin{alignat}{2}
        - \nabla \cdot (D \nabla u) & = f & \quad & \text{in } \domain{},    \\
        u                           & = g & \quad & \text{on } \domainbnd{},
    \end{alignat}
\end{subequations}
where $D := D(x) \geq d > 0$ is a scalar function.

To convert this problem into the \DD{} framework on the extended domain, $\extdomain{} \supset \domain{}$,
it is derived from extending the integrals
in the variational formulation using the characteristic function, $\chi_\domain{}$,
and surface delta function, $\delta_{\domainbnd{}}$,
\begin{align}
    \int_{\domain{}}{... \; dx}    & \mapsto \int_{\extdomain{}}{\chi_\domain{} ...\; dx}         \\
    \int_{\domainbnd{}}{... \; dS} & \mapsto \int_{\extdomain{}}{\delta_{\domainbnd{}} ...\; dx}.
\end{align}
Therefore, we need to construct a phase-field function $\phi$ satisfying
\begin{equation}
    \phi(x)
    \approx \chi_{\domain{}}
    = \begin{cases}
        1 & x \in \domain{};                      \\
        0 & x \in \extdomain{}\setminus\domain{}.
    \end{cases}
\end{equation}

\begin{defn}
    \label{defn:phi}
    We approximate the sharp characteristic function with the following smooth phase-field function:
    \begin{equation}
        \label{eq:phi}
        \phi(x) = \frac{1}{2} \left( 1 - \tanh{\left( \frac{3 r(x)}{\epsilon} \right)}\right),
        \quad \forall x \in \extdomain{}.
    \end{equation}
    Here $r(x)$ is the signed distance function for the original domain.
    Note that $\epsilon$ is a small parameter to determine the width of the smooth interfacial region.
    The above $\phi$ gives an interface width of approximately $2\epsilon$.
\end{defn}

\begin{lem}
    \label{lem:phi_prop}
    It is simple to show that the gradient of this phase-field is
    \begin{equation}
        \nabla \phi  = - \frac{6}{\epsilon} \phi(x) (1-\phi(x))\nabla r(x).
    \end{equation}
\end{lem}

Focusing on methods explored by \citet{Li2009,Yu2020}, %
these authors suggest two core ideas for constructing approximations for the Poisson problem with $D=1$ and Dirichlet boundary conditions.
\begin{subequations}
    \begin{align}
        \text{DDM1} & : & -\nabla \cdot (\phi \nabla u) & + \epsilon^{-3}(1-\phi)(u-\overline{g}) = \phi \overline{f}, \label{eq:DDM1} \\
        \text{DDM2} & : & -\phi \Delta u                & + \epsilon^{-2}(1-\phi)(u-\overline{g}) = \phi \overline{f}. \label{eq:DDM2} %
    \end{align}
\end{subequations}
This requires defining extensions $\overline{f}$ and $\overline{g}$ of $f$ and $g$, respectively.
We smoothly extend both functions so that they are constant in the normal direction to the original boundary, $\domainbnd{}$,
extending $g$ both into $\domain{}$ and $\domain{}^c$, and $f$ only into $\domain{}^c$.

Following \citet{Li2009},
the surface delta is approximated by,
\begin{equation}
    \label{eq:surfacedelta}
    \delta_{\domainbnd{}} \approx |\nabla\phi|.
\end{equation}
Note that this has the property,
$\nabla \phi = -\delta_{\domainbnd{}} n$.
This allows the outward unit normal of $\domainbnd{}$ to be approximated as
\begin{equation}
    \label{eq:normal}
    \overline{n}(x) = -\frac{\nabla\phi(x)}{|\nabla\phi(x)|}.
\end{equation}
We next define the extensions for the forcing $\overline{f}$, and boundary data $\overline{g}$:
\begin{defn}
    \label{eq:extension}
    We define the extensions for the forcing $f$ and boundary data $g$ using
    \begin{equation}
        \label{eq:gbar}
        \overline{g} = g{\left(x - r(x) \overline{n}(x) \right)},
    \end{equation}
    and
    \begin{equation}
        \label{eq:fbar}
        \overline{f} = \begin{cases}
            f(x)                                   & x \in \domain{},                         \\
            f\left(x - r(x) \overline{n}(x)\right) & x \in \extdomain{} \backslash \domain{}.
        \end{cases}
    \end{equation}
    The same is also applied to any non-constant coefficients,
    such as the diffusion coefficient $D$,
    to ensure they are extended into $B$.
\end{defn}

\begin{rmk}
    This extension is only required in a small neighbourhood around the boundary.
    The analysis by \citet{Li2009} states the neighbourhood
    ${x \in \extdomain{} : |r(x)| \ll 1}$,
    and \citet{Yu2020} states this extension must be a minimum width of $3\epsilon$.
    For simplicity, we will denote both the original functions and the extensions the same,
    unless stated otherwise.
\end{rmk}
\subsection{Weak Form Derivation}
\label{subsection:approach1}

First we will look at the derivation of DDM1 \cref{eq:DDM1},
which involves transforming the weak form.
The model Poisson problem will become:
Find $u \in H^1_0(\domain{}) + g$ such that
\begin{equation}
    \int_{\domain{}}{D \nabla u \cdot \nabla v \;dx}
    =
    \int_{\domain{}}{fv \;dx},
    \quad \forall v \in H^1_0(\domain{}).
\end{equation}
Introducing the phase-field function to extend the domain into $\extdomain{}$ gives:
Find $u \in H^1_0(\extdomain{}) + g$ such that
\begin{equation}
    \int_{\extdomain{}}{\phi D \nabla u \cdot \nabla v \;dx}
    + \int_{\extdomain{}}{{\BC{}} v \;dx}
    =
    \int_{\extdomain{}}{\phi fv \;dx},
    \quad \forall v \in H^1_0(\extdomain{}),
\end{equation}
where ${\BC{}}$ is an additional term to incorporate the boundary conditions,
and defines what PDE is solved in the external domain.
Heuristically,
the idea is to approximately solve the original PDE in $\domain{}$
and set $u=g$ in the extended region $\extdomain{}\setminus\domain{}$,
which are both combined into a single equation on $\extdomain{}$.
Using integration by parts with
the fact that $v$ vanishes on the boundary $\extdomainbnd{}$,
we see that
\begin{equation}
    \int_{\extdomain{}}{- \nabla \cdot \left(\phi D \nabla u\right) v \;dx}
    + \int_{\extdomain{}}{{\BC{}}v \;dx}
    =
    \int_{\extdomain{}}{\phi fv \;dx}.
\end{equation}
Many alternatives exist for ${\BC{}}$; however, following DDM1,
we take
\begin{equation}
    \label{eq:BC}
    {\BC{}} \propto (1-\phi)(u-g).
\end{equation}
The asymptotic analysis provided by \citet{Li2009} suggests a scaling of $\epsilon^{-3}$ should be chosen.
This is supported by the results from \citet{Franz2012},
where the scaling factor of $\epsilon^{-a}$ is used and
demonstrates $\alpha \geq 2$ is necessary.
Experimental results show a benefit of scaling this in proportion to the diffusion.
Writing this method in the bilinear form gives the first method,
\begin{meth}[DDM1]
    \label{meth:ddm1}
    Find $u \in U$ such that
    \begin{align}
        \blfa{u}{v} & :=
        \int_{\extdomain{}}{\phi D\nabla u \cdot \nabla v + \frac{D}{\epsilon^{3}}(1-\phi)u v \;dx},
        \\
        \blfl{v}    & :=
        \int_{\extdomain{}}{\phi f v + \frac{D}{\epsilon^{3}}(1-\phi)g v \;dx}.
    \end{align}
    For all $v \in V$, where
    $V = H^1_0(\extdomain{})$
    and
    $U = V + g$.
\end{meth}

We can now explore some important properties of this method.
\begin{thm}
    \label{thm:ddm1_pd}
    DDM1 (\cref{meth:ddm1}) is symmetric and coercive i.e.,
    \begin{equation}
        \blfa{u}{u} \geq C \|u\|^2
    \end{equation}
    for a constant $C>0$, with the following norm on $H^1_0(\extdomain{})$,
    \begin{equation}
        \label{eq:ddm1norm}
        \|u\|^2 :=
        \int_{\extdomain{}}{
            \phi |\nabla u|^2 + \epsilon^{-3}(1-\phi)|u|^2
            \; dx}.
    \end{equation}
\end{thm}
\begin{proof}
    We can clearly see that the bilinear form is symmetric with ${\blfa{u}{v}=\blfa{v}{u}}$.
    Furthermore, it is straightforward to see coercivity.
    \begin{equation}
        \blfa{u}{u} = \int_{\extdomain{}}{
            \phi D|\nabla u|^2 + \frac{D}{\epsilon^{3}}(1-\phi)|u|^2
            \; dx}
        \geq
        d \int_{\extdomain{}}{
            \phi |\nabla u|^2 + \epsilon^{-3}(1-\phi)|u|^2
            \; dx}
        =
        C \|u\|^2.
    \end{equation}
    Recalling from \cref{eq:poisson}, $D$ is uniformly bounded below by $d$.
    So for this case, we have $C=d$.

    It remains to show that $\|u\|$ is a norm on $H^1_0(\extdomain{})$.
    As ${\phi \in (0, 1)}$ we get ${(1-\phi)\in(0,1)}$, so
    \begin{equation}
        0 = \|u\|^2 \implies 0 = |\nabla u|^2 = |u|^2.
    \end{equation}
    So $u$ must be a identically zero on $\extdomain{}$.
    It is trivial to show the remaining norm properties,
    from the similarity of this norm (\cref{eq:ddm1norm}) to ${\|\cdot\|_{H^1_0}}$.

    Note that this modification to the norm can be seen as a weighted norm,
    with $\phi$ defined as a positive function of the distance function \citep{Cavalheiro2008}.
\end{proof}

\subsection{Strong Form Derivation}
\label{subsection:approach2}

For the derivation of DDM2,
we first start with the strong form of the problem in \cref{eq:poisson}
and multiply by $\phi$ to obtain
\begin{equation}
    -\phi\nabla \cdot (D \nabla u) = \phi f.
\end{equation}
Applying the same extension into the extended domain, the weak form becomes
\begin{equation}
    \int_{\extdomain{}}{D \nabla u \cdot \nabla (\phi v) \;dx}
    + \int_{\extdomain{}}{\BC{} \; v \;dx}
    = \int_{\extdomain{}}{\phi fv \;dx}.
\end{equation}
We will use the same choice of ${\BC{}}$ as DDM1 from \cref{eq:BC},
however this time based on the same asymptotic analysis we use $\epsilon^{-2}$.
Thus, the weak form is,
\begin{meth}[DDM2]
    \label{meth:ddm2}
    Find $u \in U$ such that
    \begin{align}
        \blfa{u}{v} & :=
        \int_{\extdomain{}}{
            D\nabla u \cdot \nabla (\phi v) + \frac{D}{\epsilon^{2}}(1-\phi)u v \;dx
        },
        \\
        \blfl{v}    & :=
        \int_{\extdomain{}}{\phi f v + \frac{D}{\epsilon^{2}}(1-\phi)g v \;dx}.
    \end{align}
    For all $v \in V$, where
    $V = H^1_0(\extdomain{})$
    and
    $U = V + g$.
\end{meth}
Exploring the same properties,
we can clearly see the bilinear form is not symmetric.
However, we are unable to comment about the coercivity of this method, but
in \cref{section:Nitsche} we will derive a new method that is coercive.

Another approach to obtain a strong form equation on the extended domain is
the Smoothed Boundary Method (SBM),
a general framework is provided by \citet{Yu2009}.

We will explore the similarities of SBM derivation to the \DD{} method by
applying this to our example problem.
First multiply \cref{eq:poisson} by $\phi$ and using the product rule we get,
\begin{equation}
    - \nabla \cdot (\phi D \nabla u)
    + \nabla\phi\cdot (D \nabla u)
    =
    \phi f.
\end{equation}
Repeating that step, results in
\begin{equation}
    - \phi \nabla \cdot (\phi D \nabla u)
    + D \left[ \nabla\phi\cdot\nabla(\phi u)
        - u|\nabla\phi|^2 \right]
    = \phi^2 f.
\end{equation}
Recall, we approximate the surface delta with $\delta_{\domainbnd{}} = |\nabla\phi|$.
Therefore, we can use the Dirichlet boundary conditions $u=g$ in that term to give
\begin{equation}
    \label{eq:SBM}
    - \phi \nabla \cdot (\phi D \nabla u)
    + D \nabla\phi \cdot \nabla(\phi u)
    - D |\nabla\phi|^2 g
    = \phi^2 f.
\end{equation}

Typically, SBM and \DD{} methods require a small offset to $\phi$ from $0$
to make sure the problem is defined and solvable on the exterior domain.
To further explicitly define the external problems,
we can add the same $\BC{}$ term in \cref{eq:BC} from DDM1/DDM2 to all methods used.

While rigorous asymptotic analysis for this modification is a topic for future work,
numerical experiments show this maintains the accuracy of the method.
Also, the experiments show that a scaling of $\epsilon^{-2}$ remains stable
and provides the lowest errors compared to other penalty parameters.
Therefore, this gives the following weak form:
\begin{meth}[SBM]
    \label{meth:SBM}
    Find $u \in U$ such that
    \begin{align}
        \blfa{u}{v} & :=
        \int_{\extdomain{}}{
            \phi D \nabla u \cdot \nabla (\phi v)
            + D \nabla\phi \cdot \nabla(\phi u) v
            + \frac{D}{\epsilon^{2}}(1-\phi) u v
            \;dx},
        \\
        \blfl{v}    & :=
        \int_{\extdomain{}}{
            \phi^2 f v
            + D |\nabla\phi|^2 g v
            + \frac{D}{\epsilon^{2}}(1-\phi) g v
            \;dx}.
    \end{align}
    For all $v \in V$, where
    $V = H^1_0(\extdomain{})$
    and
    $U = V + g$.
\end{meth}

\begin{thm}
    \label{thm:sbm_pd}
    SBM (\cref{meth:SBM}) is coercive on $H^1_0(\extdomain{})$ but not symmetric,
    with the norm
    \begin{equation}
        \label{eq:sbmnorm}
        \|u\|^2 :=
        \int_{\extdomain{}}{
            |\nabla (\phi u) |^2 + \epsilon^{-2}(1-\phi)|u|^2
            \; dx}.
    \end{equation}
\end{thm}
\begin{proof}
    The first two terms of the bilinear form clearly show it is not symmetric.

    To show this is coercive, we can simply see that when $v=u$ the bilinear is
    \begin{equation}
        \blfa{u}{u}
        =
        \int_{\extdomain{}}{
            D {|\nabla (\phi u)|}^2
            + \frac{D}{\epsilon^{2}}(1-\phi) |u|^2
            \;dx}
        \geq
        d\int_{\extdomain{}}{
            {|\nabla (\phi u)|}^2
            + \frac{1}{\epsilon^{2}}(1-\phi) |u|^2
            \;dx}
        = C \|u\|^2.
    \end{equation}

    Finally, to show \cref{eq:sbmnorm} is a norm,
    we follow from the results in (\cref{thm:ddm1_pd}).
    Given all terms are non-negative,
    \begin{equation}
        0 = \|u\|^2 \implies 0 = |\nabla (\phi u) |^2 = \epsilon^{-2}(1-\phi)|u|^2.
    \end{equation}
    Therefore, as $\epsilon^{-2}(1-\phi)>0$,
    we have $0 = |u|^2$,
    which means $u$ must be identically zero on $\extdomain{}$.
\end{proof}

\section{Diffuse Domain Nitsche's Method}
\label{section:Nitsche}

Nitsche's method is a classic method of
weakly enforcing a Dirichlet boundary condition on a PDE
\citep{Burman2012, Juntunen2009}.
A general framework is presented by \citet{Benzaken2024},
including for the Poisson model \cref{eq:poisson}.
This approach adds a penalty term instead of
enforcing the boundary condition directly in the function space.
Therefore, allowing us to apply integration by parts,
while keeping the consistency term to obtain:
\begin{equation}
    \int_{\domain{}}{D \nabla u \cdot \nabla v\;dx}
    -
    \int_{\domainbnd{}}{D \nabla u \cdot n v\;dS}
    +
    \beta\int_{\domainbnd{}}{(u-g) v \;dS}
    =
    \int_{\domain{}}{f v\;dx},
\end{equation}
where $\beta$ is a penalty coefficient used to weakly enforce the boundary condition.
Note that to make this method symmetric we can add this optional term
\begin{equation}
    -\int_{\domainbnd{}}{D \nabla v \cdot n (u-g)\;dS}.
\end{equation}
This is summarised as,
\begin{equation}
    \int_{\domain{}}{D \nabla u \cdot \nabla v\;dx}
    - \int_{\domainbnd{}}{D \nabla u \cdot n v\;dS}
    - \alpha \int_{\domainbnd{}}{D \nabla v \cdot n (u-g)\;dS}
    + \beta\int_{\domainbnd{}}{(u-g) v \;dS}
    =
    \int_{\domain{}}{fv\;dx}.
\end{equation}
The parameter $\alpha \in \{1,0,-1\}$ controls whether we have
a symmetric, non-symmetric, or skew-symmetric form.

Thus, we apply the \DD{} method using the first approach given in \cref{subsection:approach1},
with the normal
given by \cref{eq:normal}
and surface delta given in \cref{eq:surfacedelta}.
To achieve the following,
\begin{equation}
    \int_{\extdomain{}}
    \phi D \nabla u \cdot \nabla v
    + D \nabla u \cdot \nabla \phi v
    + \alpha D (u-g) \nabla \phi \cdot \nabla v
    + \beta(u-g) v |\nabla \phi|
    + \BC{}\; v
    \;dx
    =
    \int_{\extdomain{}}{\phi fv\;dx}.
\end{equation}
\begin{rmk}
    When we take $\alpha = \beta = 0$, Nitsche's method matches DDM2,
    so this method can be derived via \cref{subsection:approach2} as well.
    For Nitsche's method to be coercive, we require $\beta > 0$ in the case $\alpha \geq 0$.
    This suggests we can create a coercive DDM2 style method by
    adding these stabilising terms from Nitsche's methods in the \DD{} framework.
\end{rmk}
Based on numerical experiments we have observed that
using the same $\BC{}$ from DDM2 (\cref{meth:ddm2}) leads to the highest accuracy.
Therefore, by combining the first two terms, we use this weak form
\begin{equation}
    \int_{\extdomain{}}
    D \nabla u \cdot \nabla (\phi v)
    + \frac{D}{\epsilon^{2}}(1-\phi) (u-g) v
    + \alpha D (u-g) \nabla \phi \cdot \nabla v
    + \beta(u-g) v |\nabla \phi|
    \;dx
    =
    \int_{\extdomain{}}{\phi fv\;dx}.
\end{equation}

We will focus on the following two methods:
a stabilised method, $\alpha=0$,
and a stabilised symmetric method, $\alpha=1$,
with ${\beta=\frac{3D}{2\epsilon}\left[(\alpha + 1 )^2 (1-\phi) + (1-\phi)^2 \right]}$.
\begin{meth}[\NDDM{}]
    \label{meth:nddm}
    Find $u \in U$ such that
    \begin{align}
        \blfa{u}{v} & :=
        \int_{\extdomain{}}
        D \nabla u \cdot \nabla (\phi v)
        + \frac{D}{\epsilon^{2}}(1-\phi)u v
        +\frac{3D}{2\epsilon}\left[(1-\phi) + (1-\phi)^2\right]|\nabla\phi|uv
        \;dx,
        \\
        \blfl{v}    & :=
        \int_{\extdomain{}}
        \phi f v
        + \frac{D}{\epsilon^{2}}(1-\phi)g v
        +\frac{3D}{2\epsilon}\left[(1-\phi) + (1-\phi)^2\right]|\nabla\phi|gv
        \;dx.
    \end{align}
    For all $v \in V$, where
    $V = H^1_0(\extdomain{})$
    and
    $U = V + g$.
\end{meth}

\begin{meth}[\NSDDM{}]
    \label{meth:nsddm}
    Find $u \in U$ such that
    \begin{align}
        \blfa{u}{v} & :=
        \int_{\extdomain{}}
        D \nabla u \cdot \nabla (\phi v)
        + \frac{D}{\epsilon^{2}}(1-\phi)u v
        + D u \nabla \phi \cdot \nabla v
        +\frac{3D}{2\epsilon}\left[4(1-\phi) + (1-\phi)^2\right]|\nabla\phi|uv
        \;dx,
        \\
        \blfl{v}    & :=
        \int_{\extdomain{}}
        \phi f v
        + \frac{D}{\epsilon^{2}}(1-\phi)g v
        + D g \nabla \phi \cdot \nabla v
        +\frac{3D}{2\epsilon}\left[4(1-\phi) + (1-\phi)^2\right]|\nabla\phi|gv
        \;dx.
    \end{align}
    For all $v \in V$, where
    $V = H^1_0(\extdomain{})$
    and
    $U = V + g$.
\end{meth}

\begin{thm}
    \label{thm:nddm}
    \NDDM{} (\cref{meth:nddm}) and  \NSDDM{} (\cref{meth:nsddm}) are coercive,
    with respect to the norm
    \begin{equation}
        \|u\|^2 :=
        \int_{\extdomain{}}{
            \phi|\nabla u |^2 + \epsilon^{-2}(1-\phi)|u|^2
            \; dx}.
    \end{equation}
\end{thm}

\begin{proof}
    First, taking $v=u$ gives
    \begin{equation}
        \blfa{u}{u} \geq
        \int_{\extdomain{}}
        D \left(
            \phi |\nabla u|^2
            + \epsilon^{-2}(1-\phi)|u|^2
            + (\alpha + 1) u \nabla u \cdot \nabla \phi
            +\frac{3}{2\epsilon}\left[(\alpha + 1 )^2 (1-\phi) + (1-\phi)^2 \right]|\nabla\phi||u|^2
        \right)
        \;dx.
        \label{eq:nsddmproof1}
    \end{equation}
    Recalling from \cref{lem:phi_prop}
    the equation for $\nabla \phi$ and that $|\nabla r| = 1$,
    we replace a $|\nabla\phi|$ with $\frac{6}{\epsilon} (1-\phi)\phi$.
    Then using Young's inequality on the third term,
    with weight $w>0$,
    to get the following:
    \begin{align}
        \int_{\extdomain{}}
        D (\alpha + 1) u \nabla u \cdot \nabla \phi
        \; dx
         & \geq
        - \int_{\extdomain{}}
        D (\alpha + 1) |u| |\nabla u| |\nabla \phi|
        \; dx
        \\
         & \geq
        - \int_{\extdomain{}}
        D (\alpha + 1) |u| |\nabla u| \frac{6}{\epsilon} (1-\phi)\phi
        \; dx
        \\
         & \geq
        - \int_{\extdomain{}}
        D \left(\frac{\alpha + 1}{\sqrt{w}} \sqrt{\phi} |\nabla u| \right)
        \left(\frac{6\sqrt{w}}{\epsilon} (1-\phi)\sqrt{\phi}|u| \right)
        \; dx
        \\
         & \geq
        - \int_{\extdomain{}}
        D \frac{(\alpha + 1)^2}{2w} \phi |\nabla u|^2 +
        D \frac{18w}{\epsilon^2} (1-\phi)^2 \phi |u|^2
        \;dx.
    \end{align}

    Thus, substituting this into \cref{eq:nsddmproof1} and simplifying, we arrive at:
    \begin{align}
        \blfa{u}{u}
         & \geq
        \begin{multlined}
            \int_{\extdomain{}}
            D \left(
                \phi |\nabla u|^2
                + \epsilon^{-2}(1-\phi)|u|^2
                - \frac{(\alpha + 1)^2}{2w} \phi |\nabla u|^2
                - \frac{18w}{\epsilon^2} (1-\phi)^2 \phi |u|^2
                \right.
                \\
                \left.
                + \frac{9}{\epsilon^2}\left[
                    (\alpha + 1 )^2 (1-\phi) + (1-\phi)^2
                \right](1-\phi)\phi|u|^2
            \right)
            \;dx.
        \end{multlined}
        \\
         & \geq
        \begin{multlined}
            \int_{\extdomain{}}
            D \left(
                \left( 1 - \frac{(\alpha + 1)^2}{2w} \right)\phi |\nabla u|^2
                \right.
                \\
                \left.
                + \frac{1-\phi}{\epsilon^2}|u|^2
                \left[
                    1
                    + 9 (1-\phi) \phi
                    \left(
                        (\alpha + 1 )^2
                        + (1-\phi)
                        - 2w
                    \right)
                \right]
            \right)
            \;dx.
        \end{multlined}
    \end{align}

    To ensure coercivity we take $w = \frac{(\alpha+1)^2}{2} + \frac{1}{4}$,
    we get
    \begin{equation}
        k_1 := 1 - \frac{(\alpha + 1)^2}{2w}
        =
        1 - \frac{(\alpha+1)^2}{(\alpha+1)^2 + \frac{1}{2} }
        =
        \frac{\frac{1}{2}}{(\alpha+1)^2 + \frac{1}{2} }
        =
        \frac{1}{2(\alpha+1)^2 + 1 }
        \geq \frac{1}{9}
        > 0.
    \end{equation}
    Furthermore, for the second term,
    \begin{align}
        k_2 & := 1 + 9 (1-\phi) \phi \left((\alpha + 1 )^2 + (1-\phi) - 2w\right)
        \\
            & =
        1 + 9 (1-\phi) \phi \left( \frac{1}{2} -\phi \right)
        \\
            & =
        1 + 9 \left( \phi^3 - \frac{3}{2}\phi^2 + \frac{1}{2}\phi\right).
    \end{align}
    It is clear $k_2 = 1$ at $\phi = 0$ and $\phi = 1$.
    We can show this cubic equation has critical points at
    $\phi = \frac{1}{2} \pm \frac{\sqrt{3}}{6}$.
    Giving the minimum value within the range $0 < \phi < 1$,
    \begin{equation}
        k_2
        =
        1 + 9 \left( \phi^3 - \frac{3}{2}\phi^2 + \frac{1}{2}\phi\right)
        \geq
        1 - \frac{\sqrt{3}}{4}
        >
        0
    \end{equation}
    Consequently, as both terms positive, this gives
    \begin{equation}
        \blfa{u}{u}
        \geq
        d\min\{\frac19,1 - \frac{\sqrt{3}}{4}\}
        \int_{\extdomain{}}
        \phi |\nabla u|^2
        + \frac{1-\phi}{\epsilon^2}|u|^2
        \;dx
        \geq
        \frac{d}{9} \| u \|^2.
    \end{equation}
\end{proof}

\begin{rmk}
    This \NSDDM{} method has been explored by \citet{Nguyen2017}
    which shows a coercivity result,
    and \citet{Monte2022},
    using different derived values for $\beta$.
    However, our choice of $\beta$ only impacts the exterior with $(1-\phi)$ and is a direct result from coercivity.
\end{rmk}
\begin{rmk}
    In \ref{sec:fullnumerical}  a comparison of these stabilising terms
    in the Nitsche's method framework will show that,
    using \NDDM{} produces lower errors than the original DDM2 method.
    Also using the symmetric term in \NSDDM{},
    significantly reduces the $L^2$ error but exhibits a larger $H^1$ error,
    by almost an order of magnitude.
    It has the obvious property of producing a symmetric matrix,
    so we will focus on this going forward.
\end{rmk}

\section{Diffuse Domain Mixed Formulation}
\label{sec:mixed_methods}

\subsection{\MixedZ{}}

To derive our new methods, we start with the mixed formulation of
our model problem from \cref{eq:poisson}, given by,
\begin{subequations}
    \label{eq:mixedpoisson}
    \begin{alignat}{2}
        \sigma                    & = \nabla u, &       &                          \\
        - \nabla \cdot (D \sigma) & = f,        & \quad & \text{in } \domain{}.    \\
        u                         & = g,        & \quad & \text{on } \domainbnd{}.
    \end{alignat}
\end{subequations}

First, we use the weak form of the first equation and integrate by parts to get
\begin{equation}
    \int_{\domain{}}{\sigma \cdot \tau \; dx}
    =
    \int_{\domain{}}{\nabla u \cdot \tau \; dx}
    =
    - \int_{\domain{}}{u \nabla \cdot \tau \; dx}
    + \int_{\domainbnd{}}{g \tau \cdot n\; dS},
\end{equation}
for all test functions $\tau \in H(\mathrm{div}, \domain{})$.
Applying the \DD{} concepts for the weak form approach in \cref{subsection:approach1},
this becomes
\begin{equation}
    \label{eq:mixed0-weak1}
    \int_{\extdomain{}}{\phi \sigma \cdot \tau \; dx}
    + \int_{\extdomain{}}{\phi u \nabla \cdot \tau \; dx}
    = \int_{\extdomain{}}{\delta_{\domainbnd{}} u \tau \cdot n\; dx}
    = - \int_{\extdomain{}}{g \tau \cdot \nabla \phi\; dx},
\end{equation}
where we used \cref{eq:surfacedelta,eq:normal}.
Integrating the second term by parts on $\extdomain{}$ gives
\begin{equation}
    \int_{\extdomain{}}{\phi u \nabla \cdot \tau \; dx}
    = - \int_{\extdomain{}}{\nabla(\phi u) \cdot \tau \; dx}
    + \int_{\extdomainbnd{}}{\phi u \tau \cdot n_{\extdomain{}} \; dS}.
\end{equation}
Taking test functions $\tau \in H_0(\mathrm{div}, \extdomain{})$,
the boundary term vanishes on $\extdomainbnd{}$, giving
\begin{equation}
    \int_{\extdomain{}}{\phi \sigma \cdot \tau \; dx}
    = \int_{\extdomain{}}{\nabla(\phi u) \cdot \tau \; dx}
    - \int_{\extdomain{}}{g  \nabla \phi \cdot \tau \; dx}.
\end{equation}
Therefore, we have the following \DD{} approximation for the gradient of $u$:
\begin{equation}
    \label{eq:mixed0-strong1}
    \phi \sigma
    = \nabla(\phi u)
    - g  \nabla \phi.
\end{equation}

Moving to the second equation,
we will use an approximation based on the strong form approach in \cref{subsection:approach2},
since it has a higher $L^2$ convergence rate than the weak form approach in \cref{subsection:approach1}.
Therefore, we start with
\begin{equation}
    \label{eq:mixed0-strong2}
    -\phi \nabla \cdot (D \sigma)
    + \BC{}
    = \phi f,
\end{equation}
and applying the product rule, we obtain
\begin{equation}
    - \nabla \cdot (\phi D \sigma)
    + D \sigma \cdot \nabla \phi
    + \BC{}
    = \phi f.
\end{equation}
Multiplying by $\phi$
and substituting \cref{eq:mixed0-strong1} gives a new method,
which we will refer to as \MixedZ{}:
\begin{equation}
    \label{eq:mixed0-strongFull}
    - \phi \nabla \cdot D \left[ \nabla(\phi u ) - g \nabla \phi \right]
    + D \left[\nabla(\phi u ) - g \nabla \phi \right] \cdot \nabla \phi
    + \tBC{}
    = \phi^2 f.
\end{equation}
From experiments,
${\tBC{} = \frac{D}{\epsilon^{2}}(1-\phi)^2 (u-g)}$
was shown to be the penalty term resulting in the lowest $L^2$ and $H^1$ errors.

To make this method coercive, we add
an additional penalty term to enforce the boundary condition,
\begin{equation}
    \label{eq:mixed0-stabterm}
    \int_{\extdomain{}}{\frac{D}{4}|\nabla \phi|^2 (u-g)v\; dx}.
\end{equation}
This leads to the new method,
\begin{meth}[\MixedZ{}]
    \label{meth:mixed0}
    Find $u \in U$ such that
    \begin{align}
        \blfa{u}{v} & :=
        \int_{\extdomain{}}{
            D \nabla(\phi u ) \cdot \left[ \nabla \left(\phi v \right) + v \nabla \phi \right]
            +
            \frac{D}{\epsilon^{2}}(1-\phi)^2 u v
            + \frac{D}{4}|\nabla \phi|^2 uv
            \; dx},
        \\
        \blfl{v}    & := \int_{\extdomain{}}{
            D g \nabla \phi \cdot \left[ \nabla \left(\phi v \right) + v \nabla \phi \right]
            +\frac{D}{\epsilon^{2}}(1-\phi)^2gv
            + \frac{D}{4}|\nabla \phi|^2 gv
            + \phi^2 f v
            \; dx}.
    \end{align}
    For all $v \in V$, where
    $V = H^1_0(\extdomain)$
    and
    $U = V + g$.
\end{meth}

For \MixedZ{}, we can clearly see that this method is not symmetric.
\begin{thm}
    \label{thm:mixed0}
    \MixedZ{} (\cref{meth:mixed0}) is coercive,
    with the norm
    \begin{equation}
        \label{eq:h10normMixedZ}
        \| u \|^2 :=
        \int_{\extdomain{}} \phi^2 |\nabla u|^2 + \epsilon^{-2}(1-\phi)^2 |u|^2 \; dx.
    \end{equation}
\end{thm}
\begin{proof}
    The bilinear form with $v=u$ is
    \begin{equation}
        \blfa{u}{u}
        =
        \int_{\extdomain{}}
        D \nabla(\phi u ) \cdot [ \nabla (\phi u ) + u \nabla \phi ]
        + \frac{D}{\epsilon^2}(1-\phi)^2 |u|^2
        + \frac{D}{4}|\nabla \phi|^2 |u|^2
        \; dx.
    \end{equation}
    Using \cref{lem:phi_prop} we can split and replace a half of $(1-\phi)^2$,
    \begin{equation}
        |\nabla \phi|
        = \frac{6}{\epsilon} (1-\phi)\phi
        \leq \frac{6}{\epsilon} (1-\phi)
        \implies
        \frac{1}{36} |\nabla \phi|^2 \leq \frac{1}{\epsilon^2} (1-\phi)^2
    \end{equation}
    This gives,
    \begin{equation}
        \blfa{u}{u}
        \geq
        \int_{\extdomain{}}
        D \nabla(\phi u ) \cdot [ \nabla (\phi u ) + u \nabla \phi ]
        + \frac{D}{2 \epsilon^2}(1-\phi)^2 |u|^2
        + \frac{1}{72} D |\nabla \phi|^2 |u|^2
        + \frac{D}{4}|\nabla \phi|^2 |u|^2
        \; dx.
    \end{equation}

    Expanding the first term and replacing $\nabla(\phi u ) = \phi \nabla u  + u \nabla \phi$,
    \begin{align}
         & =
        \int_{\extdomain{}}
        D [\phi \nabla u  + u \nabla \phi] \cdot [ \phi \nabla u  + 2 u \nabla \phi ]
        + \frac{D}{2 \epsilon^2}(1-\phi)^2 |u|^2
        + \frac{D}{72} |\nabla \phi|^2 |u|^2
        + \frac{D}{4}|\nabla \phi|^2 |u|^2
        \; dx,
        \\
         & =
        \int_{\extdomain{}}
        D |\phi \nabla u|^2
        + 3D  \phi \nabla u \cdot u \nabla \phi
        + 2 D |u \nabla \phi|^2
        + \frac{D}{2 \epsilon^2}(1-\phi)^2 |u|^2
        + \frac{D}{72}|\nabla \phi|^2 |u|^2
        + \frac{D}{4}|\nabla \phi|^2 |u|^2
        \; dx,
        \\
         & =
        \int_{\extdomain{}}
        D \phi^2 |\nabla u|^2
        + 3 D \phi \nabla u \cdot u \nabla \phi
        + \frac{D}{2 \epsilon^2}(1-\phi)^2 |u|^2
        + \frac{163D}{72} |\nabla \phi|^2 |u|^2
        \; dx
    \end{align}

    Applying Young's inequality with parameter $w$:
    \begin{equation}
        3 D \phi \nabla u \cdot u \nabla \phi
        \ge
        - \frac{3 D}{2w} \phi^2 | \nabla u|^2
        - \frac{3 D w}{2} |u|^2 | \nabla \phi|^2.
    \end{equation}
    Therefore,
    \begin{align}
        \blfa{u}{u}
         & \geq
        \int_{\extdomain{}}
        D \phi^2 |\nabla u|^2
        - \frac{3 D}{2w} \phi^2 | \nabla u|^2
        - \frac{3 D w}{2} |u|^2 | \nabla \phi|^2
        + \frac{D}{2 \epsilon^2}(1-\phi)^2 |u|^2
        + \frac{163D}{72} |\nabla \phi|^2 |u|^2
        \; dx,
        \\
         & =
        \int_{\extdomain{}}
        D (1 - \frac{3}{2w}) \phi^2 |\nabla u|^2
        + \frac{D}{2 \epsilon^2}(1-\phi)^2 |u|^2
        + D (\frac{163}{72} - \frac{3 w}{2}) |\nabla \phi|^2 |u|^2
        \; dx.
    \end{align}
    Choosing, $w = \frac{163}{108} $ gives $\frac{163}{72} - \frac{3 w}{2} = 0$,
    and
    \begin{align}
         & =
        \int_{\extdomain{}}
        \frac{D}{163} \phi^2 |\nabla u|^2
        + \frac{D}{2 \epsilon^2}(1-\phi)^2 |u|^2
        \; dx,
        \\
         & \geq
        d \int_{\extdomain{}}
        \frac{1}{163} \phi^2 |\nabla u|^2
        + \frac{1}{2 \epsilon^2}(1-\phi)^2 |u|^2
        \; dx,
        \\
         & \geq
        C \|u\|^2.
    \end{align}

    Recalling \cref{eq:ddm1norm} for the norm used in the \DDMO{} in \cref{meth:ddm1},
    the norm required for \MixedZ{} is identical, except without the squared weighting.
    Therefore, we can repeat the same proof to show this is a norm.
\end{proof}
\begin{rmk}
    Our proof depends on the $\tanh$ profile used for $\phi$ and $|\nabla
        r|=1$. A more general proof is possible if we replace
    $\frac{D}{4}|\nabla \phi|^2 uv$ with $\mu D|\nabla \phi|^2 uv$ with
    $\mu>\frac{1}{4}$ using the estimate
    $|q_1|^2+q_1\cdot q_2+\mu|q_2|^2 \geq \frac{\mu-\frac14}{\mu+2} |q_1-q_2|^2$ with
    $q_1=\nabla(\phi u)$ and $q_2=u\nabla\phi$.
\end{rmk}

\subsection{\MixedO{}}
\label{subsection:mixed1}

In all the methods we have explored,
we enforce $u=g$ in the extended domain $\extdomain{}$
by adding a penalisation term to the method.
We can extend this idea to derive another method
where we also enforce $\sigma = \nabla g$ outside the domain.
Although $u = \overline{g}$ implies $\sigma = \nabla \overline{g}$,
due to interfacial regions and phase-field at the boundary,
this additional term acts as a stronger localised penalty near the boundary.

Recall the interior domain from \cref{eq:mixed0-weak1},
\begin{equation}
    \int_{\extdomain{}}{\phi \sigma \cdot \tau \; dx}
    + \int_{\extdomain{}}{\phi u \nabla \cdot \tau \; dx}
    = - \int_{\extdomain{}}{g_1 \tau \cdot \nabla \phi\; dx}.
\end{equation}
Then, we apply the same procedure in the outside region to achieve
\begin{equation}
    \int_{\extdomain{}}{(1-\phi) \sigma \cdot \tau \; dx}
    + \int_{\extdomain{}}{(1-\phi) g \nabla \cdot \tau \; dx}
    = \int_{\extdomain{}}{g_2 \tau \cdot \nabla \phi\; dx},
\end{equation}
where the functions $g_1, g_2$ will be chosen later.
Summing these two equations gives,
\begin{equation}
    \label{eq:mixed1-weak1}
    \int_{\extdomain{}}{\sigma \cdot \tau \; dx}
    + \int_{\extdomain{}}{\left[\phi u + (1-\phi) g \right] \nabla \cdot \tau \; dx}
    = \int_{\extdomain{}}{(g_2- g_1) \tau \cdot \nabla \phi\; dx},
\end{equation}
Using the same test space $\tau \in H_0(\mathrm{div}, \extdomain{})$,
integrating by parts again eliminates the boundary term,
giving the strong form
\begin{equation}
    \label{eq:mixed1-strong1}
    \sigma
    =  \nabla \left[\phi u + (1-\phi) g \right]
    + (g_2- g_1) \nabla \phi.
\end{equation}

Unlike previously in \MixedZ{}, we obtain an equation for $\sigma$ alone;
this allows us to directly substitute $\sigma$ into the \DD{} formulation of
the second equation in \cref{eq:mixed0-strong2},
\begin{equation}
    \label{eq:mixed1-strongFull}
    -\phi \nabla \cdot D \left(
        \nabla \left[\phi u + (1-\phi) g \right]
        + (g_2- g_1) \nabla \phi
    \right)
    + \frac{D}{\epsilon^{2}}(1-\phi)(u-g)
    = \phi f.
\end{equation}

Note that we have used arbitrary functions $g_1$ and $g_2$, because
we can use different assumptions on the boundary.
We could take $g_1=g_2=g$, since we know the Dirichlet boundary conditions match for both regions in the domain.
This gives the weak form of \cref{eq:mixed1-strongFull} as
\begin{equation}
    \label{eq:mixed1short-weak}
    \int_{\extdomain{}}{
        D \nabla \left[\phi u + (1-\phi) g \right] \cdot  \nabla(\phi v)
        + \frac{D}{\epsilon^{2}}(1-\phi)(u-g) v
        \; dx}
    =
    \int_{\extdomain{}}{
        \phi f v
        \; dx}.
\end{equation}

For our second method, we will take $g_1=g$ and $g_2=u$.
This can be interpreted as:
the interior PDE of $u$ needs to match the Dirichlet boundary $g$,
while the exterior PDE of $g$ needs to match $u$ on the boundary.
From experiments, this choice of $g_1$ and $g_2$ demonstrates smaller $L^2$ and $H^1$ error,
so we will refer to this new method as \MixedO{}.

We include the identical stabilising term from \MixedZ{} in \cref{eq:mixed0-stabterm} to ensure coercivity,
which gives the weak form of \cref{eq:mixed1-strongFull} as
\begin{meth}[\MixedO{}]
    \label{meth:mixed1}
    Find $u \in U$ such that
    \begin{align}
        \blfa{u}{v} & := \int_{\extdomain{}}
        D \left[\nabla (\phi u)  + u \nabla \phi\right] \cdot  \nabla(\phi v)
        + \frac{D}{\epsilon^{2}}(1-\phi)u v
        + \frac{D}{4}|\nabla \phi|^2 uv
        \; dx,
        \\
        \blfl{v}
                    & :=
        \int_{\extdomain{}}
        D \left[{-\nabla \left[(1-\phi) g \right]  + g \nabla \phi}\right] \cdot  \nabla(\phi v)
        +\frac{D}{\epsilon^{2}}(1-\phi)gv
        + \frac{D}{4}|\nabla \phi|^2 gv
        + \phi f v
        \; dx.
    \end{align}
    For all $v \in V$, where
    $V = H^1_0(\extdomain)$
    and
    $U = V + g$.
\end{meth}
We can clearly see that this is not symmetric.
The only difference between the two bilinear forms when $v=u$ is
\begin{equation}
    \blfa{u}{u} =
    \int_{\extdomain}
    D [ \nabla (\phi u ) + u \nabla \phi ] \cdot \nabla(\phi u )
    + \frac{D}{\epsilon^{2}} (1-\phi) |u|^2
    + \frac{D}{4}|\nabla \phi|^2 |u|^2
    \; dx
\end{equation}
between \MixedO{} and \MixedZ{}
is the scaling of the penalty term with $(1-\phi)$.

Using $1-\phi\geq (1-\phi)^2$, this stabilising term does not impact the previous proof for coercivity,
with respect to the norm
\begin{equation}
    \label{eq:h10normMixedO}
    \| u \|^2 :=
    \int_{\extdomain{}} \phi^2 |\nabla u|^2 + \epsilon^{-2}(1-\phi) |u|^2 \; dx.
\end{equation}

\section{Advection-Diffusion Problems}
\label{section:advection}

\subsection{DDM1 and DDM2}

As our goal is to apply the \DD{} approach to the Navier-Stokes equations,
we will use a splitting method, which requires solving
the semi-implicit momentum equation as an intermediate step.
Therefore, we will first solve advection-diffusion problems.
In this section, we will focus on the following extension of
the Poisson's equation model \cref{eq:poisson} that we have considered so far:
\begin{subequations}
    \label{eq:advdiff}
    \begin{alignat}{2}
        - \nabla \cdot (D \nabla u) + b \cdot \nabla u + m u & = f, & \quad & \text{in } \domain{}.    \\
        u                                                    & = g  & \quad & \text{on } \domainbnd{},
    \end{alignat}
\end{subequations}
We have added both a mass term, $m$, and an advection term with a given velocity vector, $b$.
Assuming that this is still an elliptic PDE, we require
\begin{equation} \label{eq:advCoercive}
    m - \frac{1}{2}\nabla \cdot b \geq 0
\end{equation}
leading to a coercive model.

As suggested by \citet{Li2009}, we use the following \DD{} approximation for the advection term:
\begin{equation}
    \label{eq:advLi}
    \int_{\extdomain{}}{\phi b \cdot \nabla u v\;dx}.
\end{equation}
This is a straightforward application of the approach described in \cref{subsection:approach2}
for DDM2 (\cref{meth:ddm2}).
However, this approximation of the advection term does not
guarantee that the resulting bilinear form is coercive.
\begin{exmp}
    \label{exmp:ddm1adv}
    For simplicity, we consider a 1D example of our model problem
    \begin{alignat*}{2}
        - 10^{-4} \frac{d^2 u(x)}{d x^2} + \frac{d u}{dx} & = 0, & \quad & x \in [-0.5, 0.5].
    \end{alignat*}
    We apply DDM1 for the diffusion and add
    \cref{eq:advLi} to approximate the advection term on the extended domain $B=[-0.675, 0.675]$.
    This gives the bilinear form,
    \begin{equation}
        \blfa{u}{v} :=
        \int_{-0.675}^{0.675}
        10^{-4} \phi \frac{d u}{dx} \frac{d v}{dx}
        + \frac{10^{-4}}{\epsilon^{3}}(1-\phi)u v \;dx
        + \phi \frac{d u}{dx} v
        \;dx.
    \end{equation}
    The signed distance function is given by ${r(x) = |x| - 0.5}$.
    Furthermore,
    we take ${\epsilon=0.035}$ and
    a uniform grid spacing of ${h=0.01}$,
    resulting in $7$ elements across the interface.

    The solution is approximated with a linear Lagrange finite element space with
    the basis $\{\eta^i_h\}$.
    This means we can construct a
    vector $\underline{u} = (u_1, \dots, u_N)$ representing coefficients of
    the discrete solution $u_h = \sum^N_{i=1} u_i \eta_i$,
    such that for the system matrix $A$,
    where $A_{ij} = a(\eta_i, \eta_j)$,
    we have $\underline{u}^T A \underline{u} = 0$.

    We fix an element inside
    the left boundary interface at ${x = -0.515}$ to $1$,
    and select another point in the right interface at ${x = 0.515}$,
    which can use a root-finding algorithm to find the corresponding $x$ below:
    \begin{equation}
        u_i = \begin{cases}
            1             & \text{if } i=16 ,                            \\
            0.643779\dots & \text{if } i=119 ,                           \\
            0             & \text{for all other } i \in \{1, \dots, N\}.
        \end{cases}
    \end{equation}
    This demonstrates that the system matrix $A$ is not positive definite
    and confirms that in general the bilinear form can be none coercive.

\end{exmp}

To guarantee stability of the method,
we need to add a stabilising term that
takes the difference of the inflow and outflow boundaries into account.
Focusing first on the advection and mass terms,
we obtain
\begin{equation}
    \begin{split}
        \int_{\domain{}}{
            b \cdot \nabla u v + m u v
            \;dx}
         & {}=
        \int_{\domain{}}{
            \nabla \cdot (b u) v - (\nabla \cdot b) uv + m u v
            \;dx}
        \\
         & {}=
        \int_{\domain{}}{
            - u b \cdot \nabla v + (m - \nabla \cdot b ) u v
            \;dx}
        + \int_{\domainbnd{}}{ b \cdot n uv  \;dS}
        \\
         & {}=
        \int_{\domain{}}{
            - u b \cdot \nabla v + (m - \nabla \cdot b ) u v
            \;dx}
        + \int_{\domainbnd{}_\text{out}}{ b \cdot n uv  \;dS}
        + \int_{\domainbnd{}_\text{in}}{ b \cdot n uv  \;dS}.
    \end{split}
\end{equation}
Note that ${\domainbnd{}_\text{in}}$ is the subset of $\domainbnd{}$
such that $b \cdot n < 0$,
while ${\domainbnd{}_\text{out}}$ is the subset with $b \cdot n \geq 0$.
Motivated by a Nitsche's method approach for advection \citep{Stoter2021, Bazilevs2007},
we replace $u$ in the final term with $g$.
This can be written as
\begin{equation}
    \int_{\domain{}}{
        b \cdot \nabla u v + m u v
        \;dx}
    =
    \int_{\domain{}}{
        - u b \cdot \nabla v + (m - \nabla \cdot b ) u v
        \;dx} \\
    + \int_{\domainbnd{}}{ b \cdot n uv  \;dS}
    - \int_{\domainbnd{}}{ [b \cdot n]^- (u-g)v \;dS}.
\end{equation}
Therefore, the advection term becomes the following formulation:
\begin{equation}
    \int_{\domain{}}{
        b \cdot \nabla u v + m u v
        \;dx}
    + \int_{\domainbnd{}}{ [-b \cdot n]^+ (u-g)v \;dS},
    \label{eq:advNitsche}
\end{equation}
where we have used the notation
\begin{equation}
    [w]^{+}
    =
    \begin{cases}
        w & w > 0   \\
        0 & w < 0~,
    \end{cases}\qquad\text{and}\qquad
    [w]^{-}
    =
    \begin{cases}
        w & w < 0   \\
        0 & w > 0~.
    \end{cases}\qquad
\end{equation}
Recall that in the \DD{} formulation,
we can approximate the normal $n$ by $-\nabla\phi/|\nabla\phi|$.
This gives the final integral:
\begin{equation}
    \int_{\domainbnd{}}{ [-b \cdot n]^+ (u-g)v \;dS} \approx
    \int_{\extdomain{}}{ [b \cdot \nabla(\phi)]^+ (u-g)v \;dx}~.
\end{equation}
Therefore, a \DD{} approximation to \cref{eq:advNitsche} is given by
\begin{equation}
    \int_{\extdomain{}}{
    \phi b \cdot \nabla u v
    + \phi m  u v
    + [ b \cdot \nabla\phi]^+ (u-g)v
    \;dx}.
\end{equation}
This is identical to \cref{eq:advLi} but with an additional stabilizing term
enforcing the boundary conditions on the inflow boundary.
From \cref{lem:phi_prop},
we know that the sign of $\nabla \phi$ when going from
the extended domain into the physical domain is positive,
and from the interior to the exterior is negative.
Therefore,
we can see that this is a similar method to the upwind discontinuous Galerkin method \citep{Brezzi2004},
replacing $u$ with $g$ on the inflow boundary.

We thus arrive at the following approximations of the advection-diffusion problem
when using DDM1 for the diffusion term:
\begin{meth}[DDM1 with stabilised advection]
    \label{meth:ddm1adv}
    Find $u \in U$ such that
    \begin{align}
        \blfa{u}{v} & :=
        \int_{\extdomain{}}{
        \phi D \nabla u \cdot \nabla v
        + \frac{D}{\epsilon^{3}}(1-\phi)u v
        + \phi b \cdot \nabla u v
        + {[b \cdot \nabla \phi]}^{+} u v
        + \phi m u v
        \;dx,
        }
        \\
        \blfl{v}    & :=
        \int_{\extdomain{}}{
        \phi f v
        + \frac{D}{\epsilon^{3}}(1-\phi)g v
        + {[b \cdot \nabla \phi]}^{+} g v
        \;dx}.
    \end{align}
    For all $v \in V$, where
    $V = H^1_0(\extdomain)$
    and
    $U = V + g$.
\end{meth}

\begin{thm}
    \label{thm:ddm1adv}
    DDM1 with stabilised advection remains coercive on $H^1_0(\extdomain())$
    with respect to the norm defined in \cref{thm:ddm1_pd}.
\end{thm}
\begin{proof}
    As we have already shown coerciveness in \cref{thm:ddm1_pd} for the
    diffusion term, we can just focus on the advection and mass terms:
    \begin{equation}
        \label{eq:ddmALLadv}
        \bar{a}(u,v) =
        \int_{\extdomain{}}{
        \phi b \cdot \nabla u v
        + {[b \cdot \nabla \phi]}^{+} uv
        + \phi m uv
        \;dx}
    \end{equation}
    for which we need to show that $\bar{a}(u,u)\geq 0$ under the conditions in \cref{eq:advCoercive} on the data.

    We can use integration by parts on the first term with $u=v$ and using
    that $u$ vanishes on the boundary to give
    \begin{equation}
        \int_{\extdomain{}}{\phi b \cdot \nabla u u\;dx}
        =
        - \int_{\extdomain{}}{\nabla \cdot (\phi b u)  u \;dx}
        =
        \int_{\extdomain{}}{-\nabla \cdot (\phi b ) |u|^2
            - \phi b \cdot \nabla u u \;dx}.
    \end{equation}
    Then, rearranging and expanding the right-hand side, we get
    \begin{equation}
        \int_{\extdomain{}}{\phi b \cdot \nabla u u\;dx}
        =
        \frac{1}{2} \int_{\extdomain{}}{-\nabla \cdot (\phi b ) |u|^2 \;dx}
        =
        \frac{1}{2} \int_{\extdomain{}}{-\phi |u|^2  \nabla \cdot b - |u|^2 b \cdot \nabla \phi \;dx}.
    \end{equation}
    This means that overall we obtain
    \begin{equation}
        \bar{a}(u,u) =
        \int_{\extdomain{}}{
        - \frac{1}{2} b \cdot \nabla \phi |u|^2
        + {[b \cdot \nabla \phi]}^{+} |u|^2
        + \phi  \left( m - \frac{1}{2} \nabla \cdot b \right) |u|^2
        \;dx}.
    \end{equation}
    From the requirements of the original elliptic PDE, we know that this last term is non-negative.
    We can split the first term into the inflow and outflow boundary regions,
    \begin{equation}
        \begin{split}
            \int_{\extdomain{}}{
            - \frac{1}{2} b \cdot \nabla \phi |u|^2
            + {[b \cdot \nabla \phi]}^{+} |u|^2
            \;dx}
             & {} =
            \int_{\extdomain{}}{
                \left(
                    - \frac{1}{2}
                    \left[
                        {[b \cdot \nabla \phi]}^{+}
                        + {[b \cdot \nabla \phi]}^{-}
                    \right]
                    + {[b \cdot \nabla \phi]}^{+}
                \right)|u|^2
                \;dx}
            \\
             & {}=
            \frac{1}{2} \int_{\extdomain{}}{
                \left(
                    {[b \cdot \nabla \phi]}^{+}
                    - {[b \cdot \nabla \phi]}^{-}
                \right)|u|^2
                \;dx}
            \\
             & {}=
            \frac{1}{2} \int_{\extdomain{}}{
                |b \cdot \nabla \phi||u|^2
                \;dx}
            \\
             & {}\geq 0.
        \end{split}
    \end{equation}
    As the additional terms are positive,
    the previous proof for coercivity implies that this is also coercive.
\end{proof}

\begin{rmk}
    The same argument can be applied to the advection contribution of DDM2, \NDDM{}, \NSDDM{}, so we use the same advection and mass terms in \cref{eq:ddmALLadv}.
\end{rmk}

\subsection{SBM}
Following the derivation for the SBM method in \cref{meth:SBM},
we multiply the model problem \cref{eq:advdiff} by $\phi$ twice.
Thus, this gives the advection and mass terms as
\begin{equation}
    \label{eq:advSBM}
    \int_{\extdomain{}}{\phi^2 b \cdot \nabla u v + \phi^2 m u v\;dx}.
\end{equation}
This is similar to the previous suggestion for DDM1, \cref{eq:advLi},
and requires a stabilising term to be added.

\begin{exmp}
    Repeating the previous \cref{exmp:ddm1adv} with
    SBM and the advection term in \cref{eq:advSBM}
    we have the bilinear form,
    \begin{equation}
        \blfa{u}{v} :=
        \int_{-0.675}^{0.675}
        10^{-4} \phi \frac{d u}{dx} \frac{d (\phi v) }{dx}
        + 10^{-4} \frac{d \phi}{dx} \frac{d (\phi u)}{dx} v
        \\
        + \frac{10^{-4}}{\epsilon^{2}}(1-\phi) u v
        + \phi^2 \frac{d u}{dx} v
        \;dx.
    \end{equation}
    We can construct an example for
    the coefficient vector $\underline{u} = (u_1, \dots, u_N)$
    such that $\underline{u}^T A \underline{u} = 0$.
    For example,
    \begin{equation}
        u_i = \begin{cases}
            1             & \text{if } i=16 ,                            \\
            0.900367\dots & \text{if } i=119 ,                           \\
            0             & \text{for all other } i \in \{1, \dots, N\}.
        \end{cases}
    \end{equation}

\end{exmp}

Therefore, we get the following approximation of the elliptic
problem using SBM:
\begin{meth}[SBM with stabilized advection]
    \label{meth:SBMadv}
    Find $u \in U$ such that
    \begin{align}
        \blfa{u}{v} & :=
        \int_{\extdomain{}}{
            \begin{multlined}[t][0.7\textwidth]
                \phi D \nabla u \cdot \nabla (\phi v)
                + D \nabla\phi \cdot \nabla(\phi u) v
                + \frac{D}{\epsilon^{2}}(1-\phi) u v
                + \phi^2 b \cdot \nabla u v
                \\
                + \phi{[b \cdot \nabla \phi]}^{+} u v
                + \phi^2 m u v
                \;dx,
            \end{multlined}
        }
        \\
        \blfl{v}    & :=
        \int_{\extdomain{}}{
            \phi^2 f v
            + D |\nabla\phi|^2 g v
            + \frac{D}{\epsilon^{2}}(1-\phi) g v
            + \phi{[b \cdot \nabla \phi]}^{+} g v
            \;dx,}
    \end{align}
    for all $v \in V$, where
    $V = H^1_0(\extdomain)$
    and
    $U = V + g$.
\end{meth}

\begin{thm}
    \label{thm:SBMadv}
    SBM with stabilised advection remains coercive
    with respect to the norm defined in \cref{thm:sbm_pd}.
\end{thm}
\begin{proof}
    Here we will closely follow the proof of \cref{thm:ddm1adv}.

    Coerciveness has already been shown for the diffusion term in \cref{thm:sbm_pd},
    so we will focus on the advection and mass terms:
    \begin{equation}
        \label{eq:sbmALLadv}
        \bar{a}(u,v) =
        \int_{\extdomain{}}{
            \phi^2 b \cdot \nabla u v
            + \phi {[b \cdot \nabla \phi]}^{+} uv
            + \phi^2 m u v
            \;dx}.
    \end{equation}
    We need to show that $\bar{a}(u,u)\geq 0$ under the conditions in \cref{eq:advCoercive} on the data.

    Integrating by parts the first term and rearranging with $u=v$ gives,
    \begin{equation}
        \int_{\extdomain{}}{\phi^2 b \cdot \nabla u u\;dx}
        =
        \frac{1}{2} \int_{\extdomain{}}{
            -\phi^2 |u|^2  \nabla \cdot b
            - |u|^2 b \cdot \nabla{(\phi^2)}
            \;dx}
        =
        \int_{\extdomain{}}{
            - \frac{1}{2} \phi^2 |u|^2  \nabla \cdot b
            - |u|^2 b \cdot \nabla \phi \phi
            \;dx}.
    \end{equation}
    So that overall, we obtain
    \begin{equation}
        \bar{a}(u,u) =
        \int_{\extdomain{}}{
            - \phi b \cdot \nabla \phi |u|^2
            + \phi {[b \cdot \nabla \phi]}^{+} |u|^2
            + \phi^2  \left( m - \frac{1}{2} \nabla \cdot b \right) |u|^2
            \;dx}.
    \end{equation}
    From the requirements of the original elliptic PDE, we know that this last term is non-negative.
    Focusing on the boundary terms, we can split $b \cdot \nabla \phi$ into the inflow and outflow boundary regions,
    \begin{equation}
        \begin{split}
            \int_{\extdomain{}}{
                \left(
                    - \phi b \cdot \nabla \phi
                    + \phi {[b \cdot \nabla \phi]}^{+}
                \right) |u|^2
                \;dx}
             & {} =
            \int_{\extdomain{}}{
                \left(
                    - \phi
                    \left(
                        {[b \cdot \nabla \phi]}^{+}
                        + {[b \cdot \nabla \phi]}^{-}
                    \right)
                    + \phi {[b \cdot \nabla \phi]}^{+}
                \right)|u|^2
                \;dx}
            \\
             & {}=
            \int_{\extdomain{}}{
                -\phi {[b \cdot \nabla \phi]}^{-}|u|^2
                \;dx}
            \\
             & {}=
            \int_{\extdomain{}}{
                \phi {[-b \cdot \nabla \phi]}^{+}|u|^2
                \;dx}
            \\
             & {}\geq 0.
        \end{split}
    \end{equation}
    Similarly, this shows that the stabilising term is sufficient for coercivity.
\end{proof}

\subsection{\MixedZ{}}
For the mixed methods,
we again consider the mixed formulation of our problem
\begin{subequations}
    \begin{alignat}{2}
        \label{eq:mixedpoisson-adv}
        \sigma                                        & = \nabla u, &       &                          \\
        - \nabla \cdot (D \sigma) + b\cdot\sigma + mu & = f,        & \quad & \text{in } \domain{}.    \\
        u                                             & = g         & \quad & \text{on } \domainbnd{}.
    \end{alignat}
\end{subequations}
Recall that in \MixedZ{},
we obtained the \DD{} approximation for $\sigma$ of the form
$\phi \sigma = \nabla(\phi u) -  g \nabla\phi$.
Then, we needed to
multiply the second equation by $\phi^2$.
Consequently, the \MixedZ{} approximation of the advection term is given by
\begin{equation}
    \int_{\extdomain{}}{\phi^2 b \cdot \sigma v\;dx}
    = \int_{\extdomain{}}{
        \phi b \cdot \left({\nabla(\phi u) -  g \nabla\phi}\right) v
        \;dx}.
\end{equation}
This gives the following method:
\begin{meth}[\MixedZ{} with advection]
    \label{meth:mixed0adv}
    Find $u \in U$ such that
    \begin{align}
        \blfa{u}{v} & :=
        \int_{\extdomain{}}{
            D \nabla(\phi u ) \cdot \left[ \nabla \left(\phi v \right) + v \nabla \phi \right]
            + \frac{D}{\epsilon^{2}}(1-\phi)^2 u v
            + \frac{D}{4}|\nabla \phi|^2 uv
            + \phi b \cdot \nabla(\phi u) v
            + \phi^2 m u v
            \; dx,
        }
        \\
        \blfl{v}    & := \int_{\extdomain{}}{
            D g \nabla \phi \cdot \left[ \nabla \left(\phi v \right) + v \nabla \phi \right]
            +\frac{D}{\epsilon^{2}}(1-\phi)^2gv
            + \frac{D}{4}|\nabla \phi|^2 gv
            + \phi b \cdot  g \nabla\phi v
            + \phi^2 f v
            \; dx.}
    \end{align}
    For all $v \in V$, where
    $V = H^1_0(\extdomain)$
    and
    $U = V + g$.
\end{meth}
Note that we do not have to add any special treatment for the inflow boundary,
as we can show that the terms added for the advection
are already coercive.
\begin{thm}
    \label{thm:mixed0adv}
    \MixedZ{} with advection remains coercive
    with respect to the norm defined in \cref{thm:mixed0}.
\end{thm}
\begin{proof}
    We can rearrange the advection term in the bilinear form with $v=u$:
    \begin{align}
        \int_{\extdomain{}}{
            \phi b \cdot {\nabla(\phi u)}u
            \;dx}
         & =
        \label{eq:mixedCoerciveRearrange}
        - \int_{\extdomain{}}{
            \phi^2 \nabla \cdot b |u|^2
            + \phi b \cdot \nabla (\phi u) u
            \;dx}.
    \end{align}
    Therefore, the advection and mass term become
    \begin{equation}
        \int_{\extdomain{}}{
            \phi b \cdot {\nabla(\phi u)}u
            + \phi^2 m |u|^2
            \;dx}
        =
        \int_{\extdomain{}}{
            \phi^2 \left(m - \frac{1}{2}\nabla \cdot b \right) |u|^2
            \;dx}
        \geq 0.
    \end{equation}
    This is non-negative,
    using assumption \cref{eq:advCoercive} on the data of the PDE,
    so the advection term is positive semi-definite. %
    Therefore, combining this with \cref{thm:mixed0},
    we see that \MixedZ{} with advection is stable.
\end{proof}
\begin{rmk}
    Our experiments indicate that \MixedO{} is outperformed by
    \MixedZ. Consequently we are omitting
    the proof of coercivity for the advection-diffusion problem since it
    involves the additional assumption of $m\geq 0$ or a more complicated
    form of the stabilization.
    In our practical applications this condition is naturally satisfied.
\end{rmk}

\section{Numerical Experiments}
\label{sec:Numerical}

We will now look at some numerical experiments of the new methods
to evaluate their potential.
These will be implemented in Python using
the package \ddfem{} by \citet{Benfield2025},
with meshes generated with \gmsh{} \citep{Geuzaine2009}
and \dunealu{} \citep{Alkaemper2016},
and solved using \dunefem{} \citep{Dedner2020}.
First we shall focus on
the Advection-Diffusion model problem in \cref{eq:advdiff},
using a Lagrange space on triangles.
The goal is to apply these methods to Navier-Stokes,
where we will use the Taylor-Hood space,
so we need to test both first-order and second-order basis functions.
Initially, we solve the corresponding linear system with a direct solver to
ensure errors arise only from the approximations.
This allows a manufactured solution to be used for verification,
while varying parameters will show the sources of error.

To compute the errors on the physical domain $\domain{}$, we use
\begin{subequations}
    \label{eq:errors}
    \begin{align}
        E_{L^2} & := \|u - u_h\|_{L^2(\domain)},                                          \\
        E_{H^1} & := \|u - u_h\|_{L^2(\domain)} + \|\nabla u - \sigma_h\|_{L^2(\domain)},
    \end{align}
\end{subequations}
where $\sigma_h$ is the reconstructed discrete gradient defined by \cref{eq:mixed0-strong1} for \MixedZ{}
and \cref{eq:mixed1-strong1} for \MixedO{} used in the derivation of the mixed methods,
while $\sigma_h = \nabla u_h$ for all other methods.
Recall from \cref{subsection:sharp2dd} that these integrals over $\domain{}$
are evaluated on $\extdomain{}$ via the characteristic function $\chi_{\domain{}}$.
Our numerical experiments show that evaluating the $H^1$ error using the reconstructed flux from the derivation of the mixed methods leads to significant improvements compared to using the gradient $\nabla u_h$ instead. For example, in our tests the $H^1$ errors using $\nabla u_h$ also in \MixedZ{} were only very slightly better then those reported here for \NSDDM{}.
The improvement is close to the boundary since $\sigma_h$ includes the (extended) boundary values $g$ in that region.
A comparison of the pointwise gradient errors for the diffusion problem on
the arc domain is given in \ref{sec:gradError}.
We could not find a similar flux reconstruction $\sigma_h$ for the methods not based on a mixed formulation.
We will also compare the experimental order of convergence (EOC), defined as
\begin{equation}
    \text{EOC} = \frac{\ln\left(E_k/E_{k+1}\right)}{\ln\left(h_k/h_{k+1}\right)},
\end{equation}
where $E_k$ and $E_{k+1}$ are the $L^2$ or $H^1$ errors (as defined in \cref{eq:errors})
for meshes with sizes $h_k$ and $h_{k+1}$ respectively.

The interface width is always chosen such that it will contain at least 7 elements
i.e., $2\epsilon = 7h$.
This choice connects the interface width parameter, $\epsilon$,
to the mesh size discretisation, $h$.
It asserts that the smooth interfacial region is resolved without numerical artefacts from a sharp jump.
Furthermore, this means as the mesh is refined, $h \to 0$, the interface width decreases, $\epsilon \to 0$,
keeping consistent with the asymptotic analysis of the \DD{} methods.
So the numerical convergence rates shown in this paper combine simultaneous
reduction of $h$ and $\epsilon$ but will be reported using $h$ as free
parameter.

Here we will focus on the two methods with the lowest errors,
but the full results are contained in \ref{sec:fullnumerical}.
All our numerical experiments show that
the original or stabilised DDM1 (\cref{meth:ddm1},\cref{meth:ddm1adv}) and DDM2 (\cref{meth:ddm2})
produce higher errors in the $L^2$ norm compared to the new methods.

\subsection{Advection-Diffusion Problems}
\label{subsection:ad_numerical}

Initially we take the advection term $b=(0,0)$,
to see the accuracy of the diffusion and boundary term.
Also, we take the non-constant scalar diffusion coefficient,
\begin{equation}
    - \nabla \cdot \left[ (1+|x|)\nabla u \right] = f,
\end{equation}
where $f$ is constructed such that we have the exact solution
\begin{equation}
    u(x) = \cos{\left(1.8 \pi x_1\right)} \cos{\left(2.6 \pi x_2\right)}.
\end{equation}

\begin{figure}[htbp]
    \centering
    \includegraphics{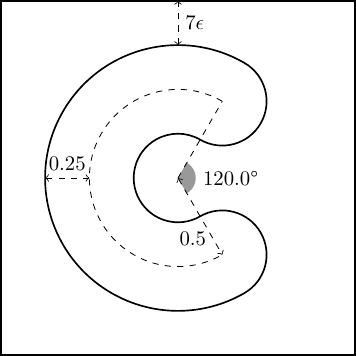}{}
    \caption{Arc shape domain.}
    \label{fig:arc}
\end{figure}

Two domains will be used:
a smooth arc shape extended into a simple square
with a width $7\epsilon$ larger than the outer radius;
and an inverted domain, which has the arc cutout from a square,
with the fixed amount extended inside the arc.
For convergence tests it is required to start with a mesh size small enough
to ensure the new external boundary conditions do not impact the results.
An example of this domain is shown in \cref{fig:arc},
and a SDF can be derived using \citet{Quilez2020}.
The use of \ddfem{} allows complex domains to easily be defined by
combining SDFs, which is demonstrated by \citet{Benfield2025}.

\begin{figure}[htbp]
    \centering
    \includegraphics{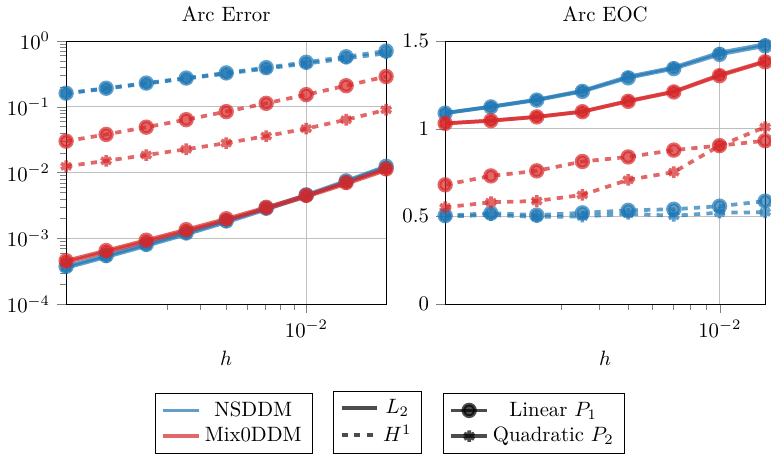}{}
    \caption{Diffusion Dominated problem error plot for Arc domain.}
    \label{fig:difArc}
\end{figure}

\begin{figure}[htbp]
    \centering
    \includegraphics{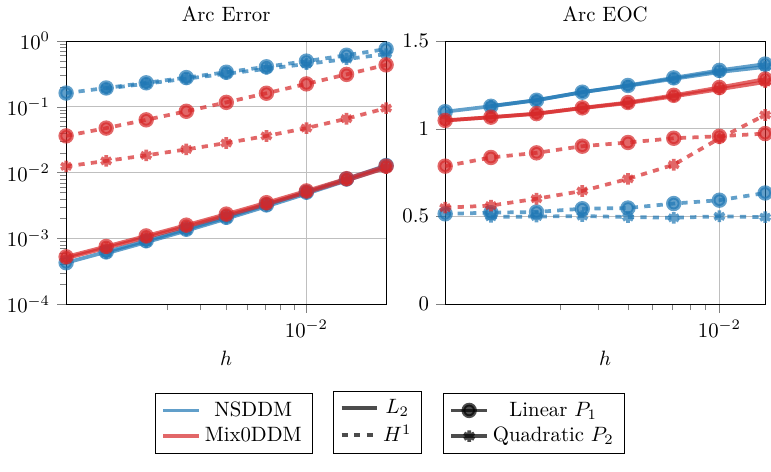}{}
    \caption{Diffusion Dominated problem error plot for Inverted Arc domain.}
    \label{fig:difInvArc}
\end{figure}

These results are shown in \cref{fig:difArc,fig:difInvArc}.
\NSDDM{} (\cref{meth:nsddm}) has the lowest $L^2$ errors except at largest $h$ values,
and is very closely followed by \MixedZ{} (\cref{meth:mixed0}).
These will be a useful combination of methods to focus on
to limit the combination of parameters,
as this explores the new mixed methods and our version using Nitsche's method.
This means we have a symmetric method that offers a high accuracy $L^2$ result,
which would be beneficial for certain problems, for example the Stokes fluid flow.

In the $H^1$ norm, \MixedZ{} produces the lowest errors of all methods tested,
and is closely followed by \MixedO{}.
However, \NSDDM{} produces larger $H^1$ errors than the other methods,
this means a compromise of the symmetric property and $H^1$ accuracy needs to be considered
depending on the problem.
Furthermore, using quadratic basis functions improves the $H^1$ errors for \MixedZ{},
whereas \NSDDM{} is not significantly affected by increasing the basis order.
We repeat the same experiments but increase the advection term significantly
to investigate the errors for an advection-dominated problem,
\begin{equation}
    - \nabla \cdot \left[ (1+|x|)\nabla u \right]
    - 100 \begin{bmatrix}
        x_2 \\ x_1
    \end{bmatrix} \cdot \nabla u = f.
\end{equation}

\begin{figure}[htbp]
    \centering
    \includegraphics{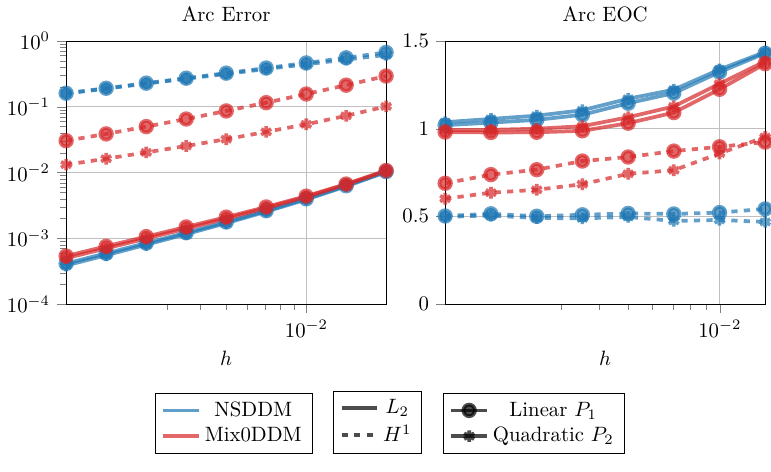}{}
    \caption{Advection Dominated problem error plot for Arc domain.}
    \label{fig:adArc}
\end{figure}

\begin{figure}[htbp]
    \centering
    \includegraphics{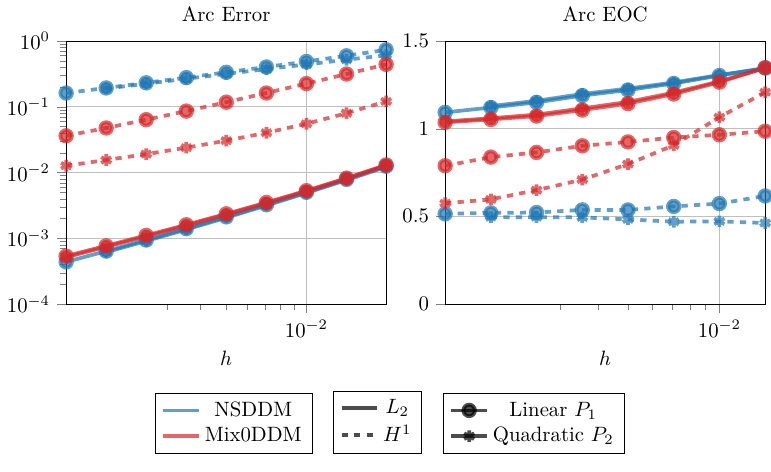}{}
    \caption{Advection Dominated problem error plot for Inverted Arc domain.}
    \label{fig:adInvArc}
\end{figure}

In \cref{fig:adArc,fig:adInvArc} we can see a very similar outcome to the case without
advection.
Again, \MixedZ{} produces only a slightly larger error in the $L^2$-norm
compared to \NSDDM{}, with no advantage gained from using a quadratic basis.
However, there is a clear difference in the $H^1$-norm
between the two methods, with \MixedZ{} showing a higher EOC than \NSDDM{}.

\subsection{Searchlight Advection Problem}

\begin{figure}[htbp]
    \centering
    \includegraphics{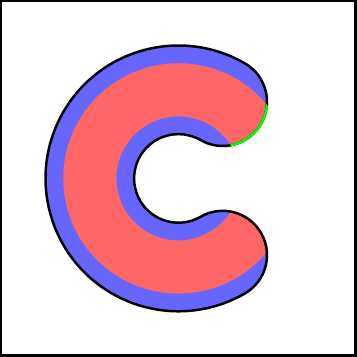}{}
    \caption{Solution for searchlight with arc domain.}
    \label{fig:arcSearchlight}
\end{figure}

To further evaluate the performance of the methods for advection dominated problems,
we will look at
\begin{subequations}
    \label{eq:searchlight}
    \begin{gather}
        - D \Delta u
        - \begin{bmatrix}
            -x_2 \\ x_1
        \end{bmatrix} \cdot \nabla u = 0.
        \\
        u = \begin{cases}
            0.5  & |x - 0.5 (\cos\frac{\pi}{6}, \sin\frac{\pi}{6} )|^2 < 0.15^2 \\
            -0.5 & \text{otherwise}
        \end{cases}
        \quad \text{on } \domainbnd{}.
    \end{gather}
\end{subequations}
The problem is shown in \cref{fig:arcSearchlight}:
the green part of the boundary shows the inflow
and black is characteristic or outflow.
While an analytical solution is not available for this problem,
this figure shows an approximation of the expected limiting solution in the
limit $D\to0$.

\begin{figure}[htbp]
    \centering
    \includegraphics{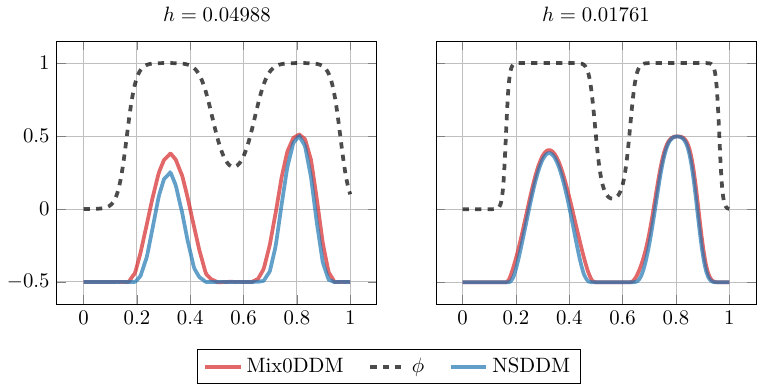}{}
    \caption{Advection problem sample for Arc domain.}
    \label{fig:sampleArc}
\end{figure}

To evaluate the results with $D=10^{-3}$,
we can sample across the yellow line in the figure,
with end points $[-0.7, -0.8]$ and $[0.1, 0.8]$.
The results are shown in \cref{fig:sampleArc} with two different spatial resolutions.
Both methods converge to the expected result
with \MixedZ{} producing values closer to the expected solution,
especially in the downwind part of the domain (left side of the plots).

\begin{figure}[htbp]
    \centering
    \begin{subfigure}[t]{\textwidth}
        \centering
        \includegraphics{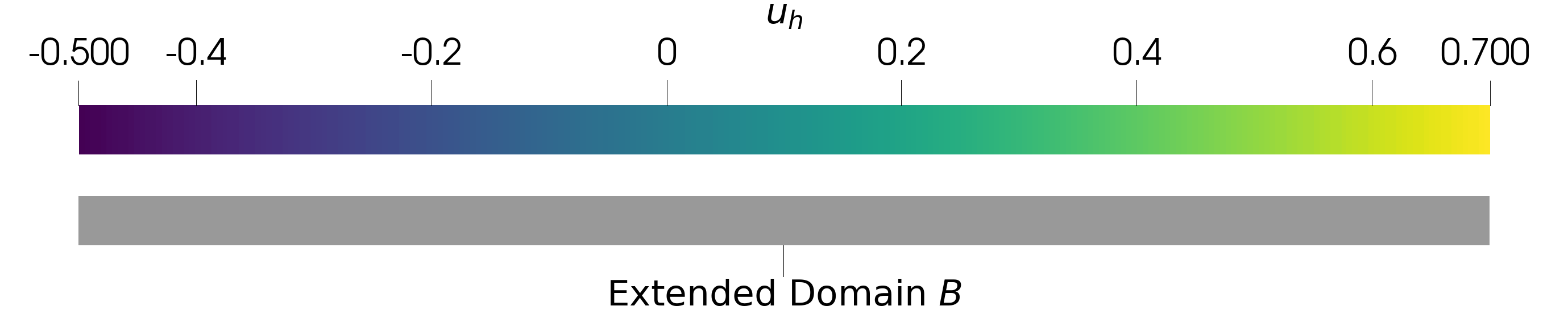}
    \end{subfigure}
    \hfill
    \begin{subfigure}{0.495\textwidth}
        \centering
        \includegraphics{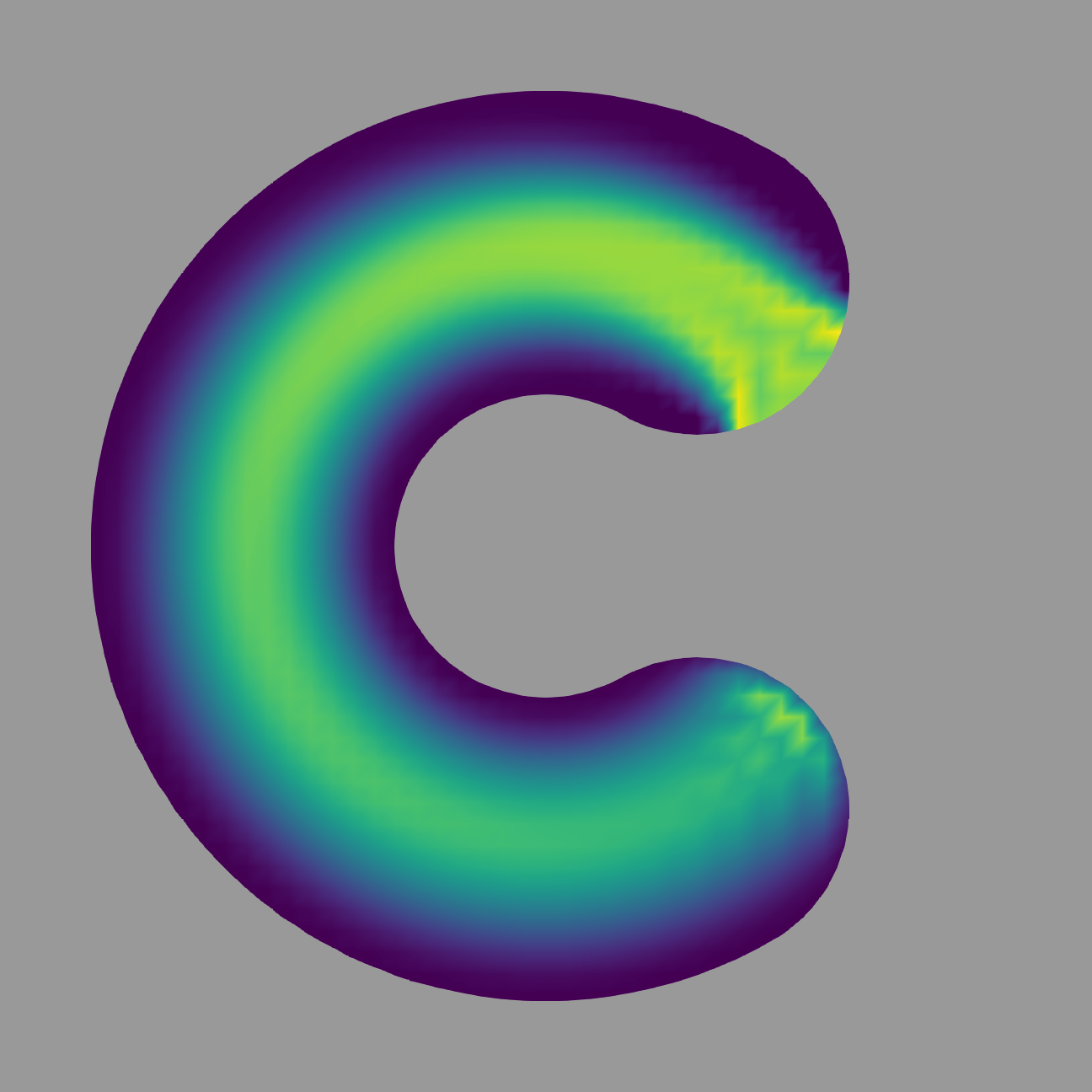}
        \caption{$h=0.04988$}
    \end{subfigure}
    \hfill
    \begin{subfigure}{0.495\textwidth}
        \centering
        \includegraphics{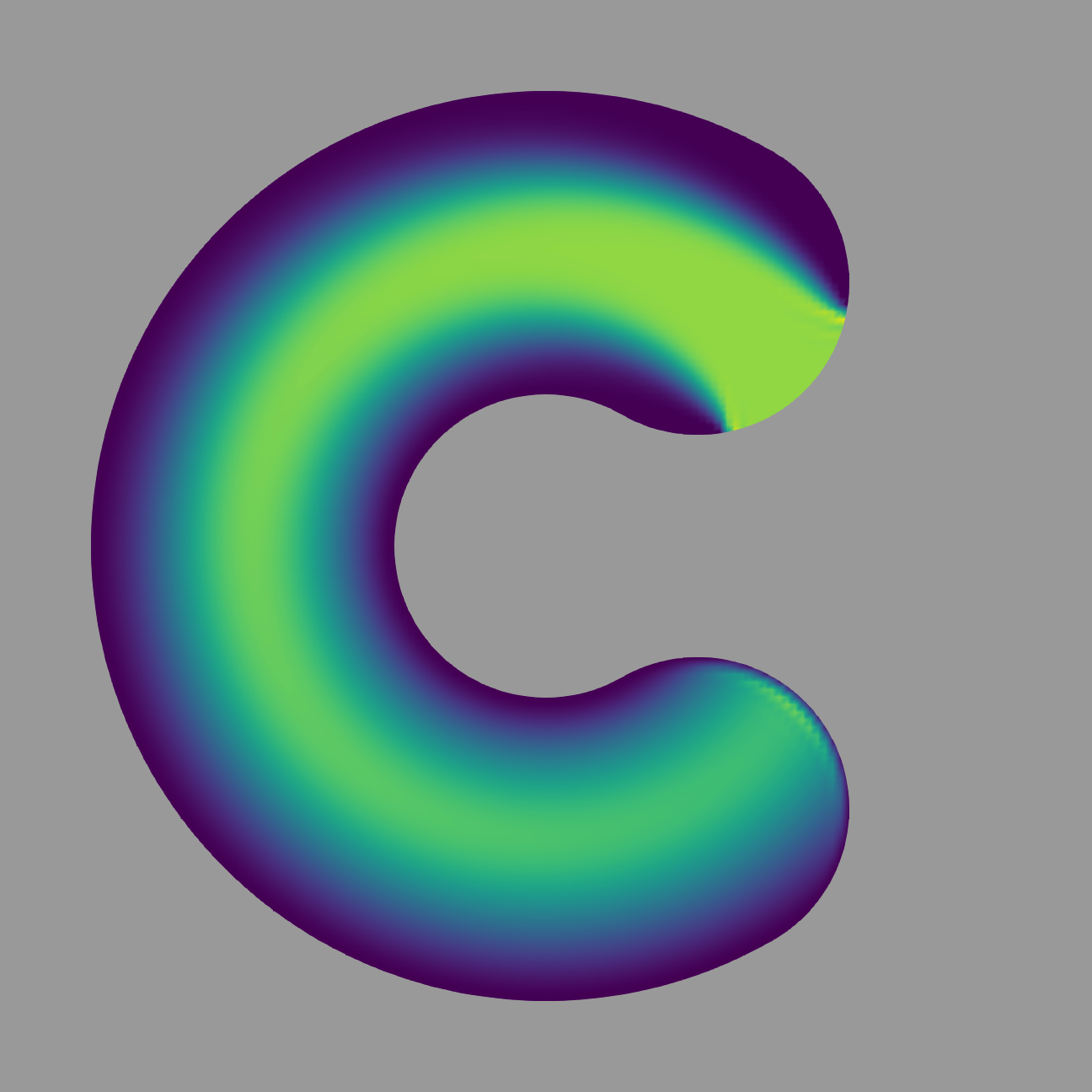}
        \caption{$h=0.01761$}
    \end{subfigure}
    \hfill
    \caption{Contour plots of $u_h \chi_{\domain{}}$ with \MixedZ{} approximation solutions for Arc domain.}
    \label{fig:Arc2dMixed0}
\end{figure}

\begin{figure}[htbp]
    \centering
    \begin{subfigure}[t]{\textwidth}
        \centering
        \includegraphics{Fig08-09Legend.png}
    \end{subfigure}
    \hfill
    \begin{subfigure}{0.495\textwidth}
        \centering
        \includegraphics{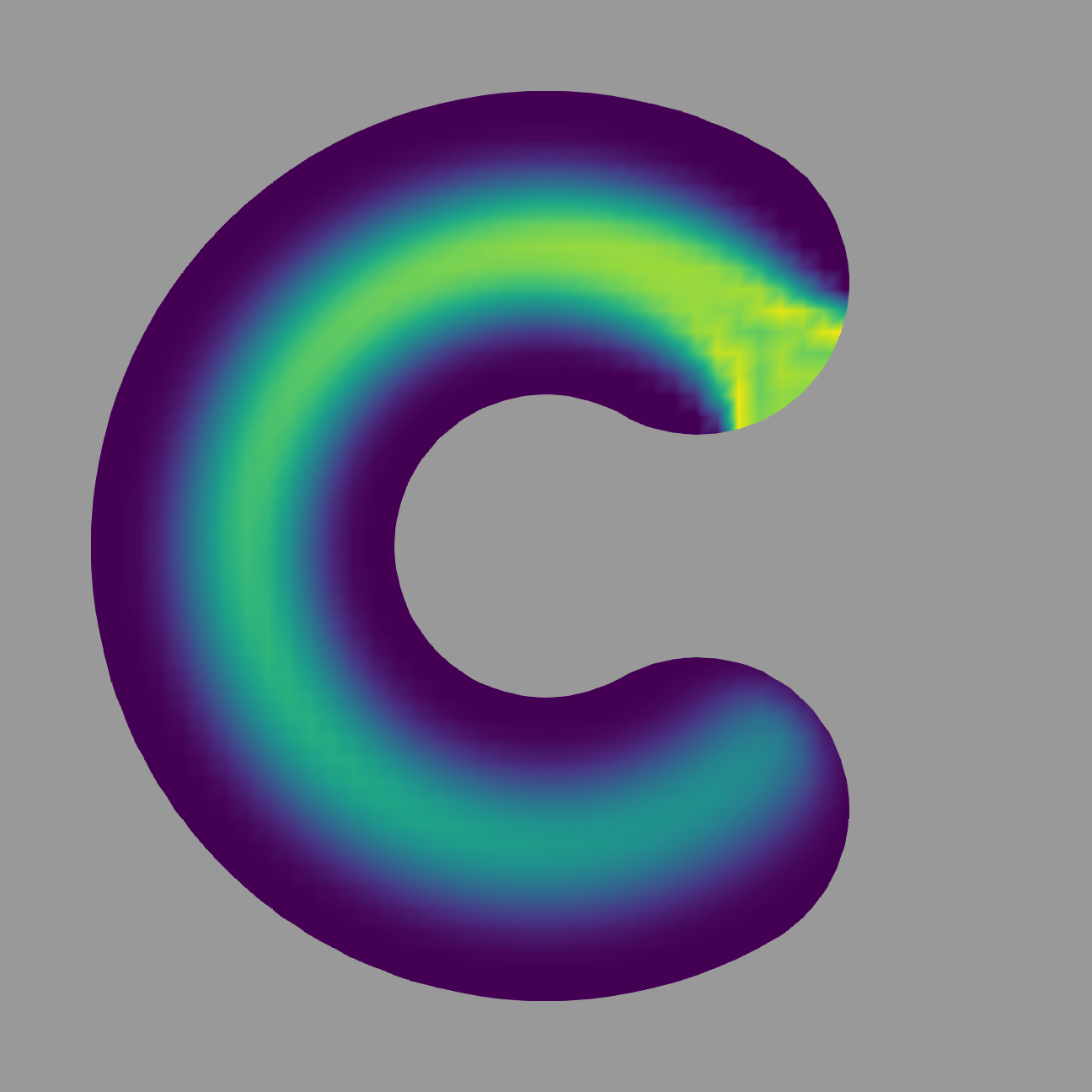}
        \caption{$h=0.04988$}
    \end{subfigure}
    \hfill
    \begin{subfigure}{0.495\textwidth}
        \centering
        \includegraphics{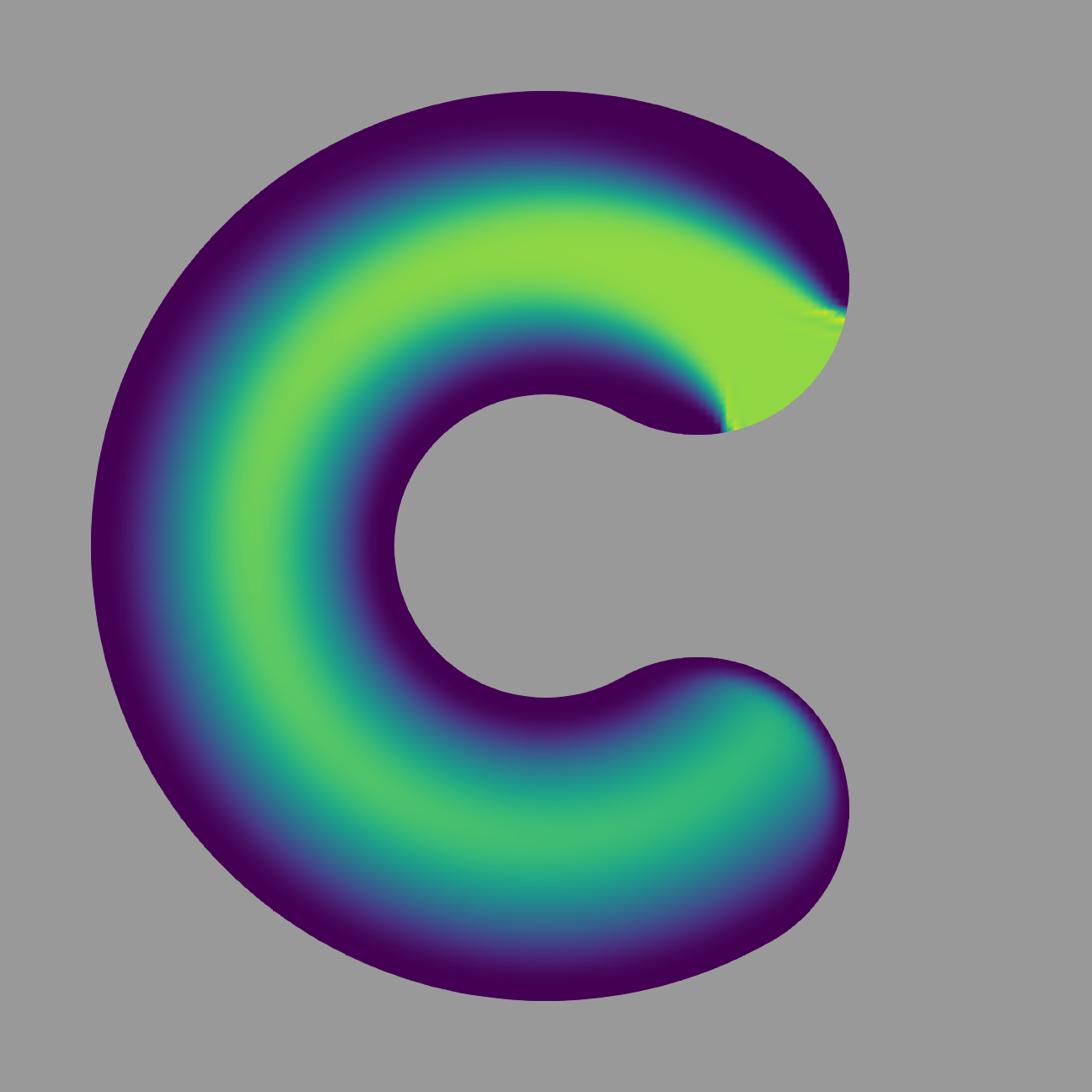}
        \caption{$h=0.01761$}
    \end{subfigure}
    \hfill
    \caption{Contour plots of $u_h \chi_{\domain{}}$ with \NSDDM{} approximation solutions for Arc domain.}
    \label{fig:Arc2dNSDDM}
\end{figure}

The results in \cref{fig:Arc2dMixed0} and \cref{fig:Arc2dNSDDM},
show a significant difference between the methods on the coarser mesh.
\MixedZ{} matches the expected solution more closely,
although it shows some oscillations at the outflow boundary.
The fine mesh shows similar results,
there is still a significant difference at the outflow boundary.
\NSDDM{} clearly enforces $u = -0.5$ more strongly but
with a larger than expected boundary layer.
\MixedZ{} seems to produce a very thin layer and
does not resolve $u = -0.5$,
more similar to what is seen in
discontinuous Galerkin methods.


\subsection{Incompressible Navier-Stokes}

Finally, we have reached our goal of the incompressible Navier-Stokes problem.
We will solve this using
incremental pressure correction (IPC),
which is a
second-order projection scheme from \citet{Guermond1998, Guermond2006}.
Note that a simpler \DD{} projection method has been used by \citet{Aland2010}.
The key benefit of this method is to decouple the velocity and pressure,
so they are solved independently.

We partition the domain boundary into
a Dirichlet boundary $\domainbnd{}_D$ (walls and inflow)
and an outflow boundary $\domainbnd{}_N$ (do-nothing).
Note, when $\domainbnd{}_N = \emptyset$,
Dirichlet conditions hold on the entire boundary $\domainbnd{} = \domainbnd{}_D$.
First, we compute a tentative velocity by removing the current pressure term in the momentum equation
and divergence-free constraint.
For a time step of $\tau$, we solve the following for $\tilde{u}^n$
\begin{subequations}
    \label{eq:HIPC1}
    \begin{alignat}{2}
        \frac{3 \tilde{u}^n - 4 u^{n-1} + u^{n-2}}{2\tau}
        - \nu \Delta \tilde{u}^n
        + ( u^{n}_{\star} \cdot \nabla)\tilde{u}^n
                               & = f^{n} - \nabla p^{n-1} & \quad & \text{in } \domain{},
        \\
        \tilde{u}^n            & = g^n                    & \quad & \text{on } \domainbnd{}_D,
        \\
        (\nabla \tilde{u}^n) n & = 0                      & \quad & \text{on } \domainbnd{}_N,
    \end{alignat}
\end{subequations}
where $u^{n}_{\star} = 2u^{n-1} - u^{n-2}$.

Then using the Helmholtz decomposition,
the tentative velocity is projected into the divergence-free space,
giving two equations to solve.
Next, solving for the pressure correction $\delta p^n$,
\begin{subequations}
    \label{eq:HIPC2}
    \begin{alignat}{2}
        - \Delta {\delta{p}}^n       & = -\frac{3}{2\tau}\nabla \cdot \tilde{u}^n & \quad & \text{in } \domain{},
        \\
        \nabla {\delta{p}}^n \cdot n & = 0                                        & \quad & \text{on } \domainbnd{}_D,
        \\
        {\delta{p}}^n                & = 0                                        & \quad & \text{on } \domainbnd{}_N.
    \end{alignat}
\end{subequations}
Then updating the pressure,
\begin{equation}
    \label{eq:HIPC2update}
    p^{n} = p^{n-1} + {\delta{p}}^n,
\end{equation}
and solving for the velocity,
\begin{equation}
    \label{eq:HIPC3}
    u^{n} = \tilde{u}^{n} - \frac{2\tau}{3} \nabla {\delta{p}}^n.
\end{equation}
\begin{rmk}
    If $\domainbnd{}_N = \emptyset$,
    homogeneous Neumann conditions hold on the entire boundary $\domainbnd{}$,
    which means we can only determine the pressure up to a constant.
\end{rmk}
Note due to the uncoupling,
this introduces a splitting error of order $\mathcal{O}(\tau^2)$,
so there is no improvement to using a higher order discretisation in time \citep[Section 7.5]{John2016}.

To solve this problem in the \DD{} approach
we use the components previously introduced.
Although, due to the advection term being divergence-free,
in \ddfem{} the advection has been implemented as a source term,
so we need to incorporate new terms for gradients and divergence.
For the gradient, we can follow the same approach as the mass term,
\begin{equation}
    \int_{\domain{}}{\nabla p \cdot v \; dx}
    \approx \int_{\extdomain{}}{ \phi \nabla p \cdot v \;dx}.
\end{equation}
Recall, \MixedZ{} (\cref{thm:mixed0adv}) and SBM (\cref{thm:SBMadv})
both require multiplying by $\phi^2$.

\Cref{eq:HIPC2} requires a Neumann boundary condition.
Neumann boundary conditions in the \DD{} framework were introduced by \citet{Li2009,Lervaag2015}.
For example, changing the boundary conditions on the model problem \cref{eq:poisson},
\begin{subequations}
    \label{eq:poisson_Neumann}
    \begin{alignat}{2}
        - \nabla \cdot (D \nabla p) & = f   & \quad & \text{in } \domain{},
        \\
        D \nabla p  \cdot n         & = p_N & \quad & \text{on } \domainbnd{},
    \end{alignat}
\end{subequations}
the \DD{} problem becomes,
\begin{equation}
    \label{eq:DDNeumann}
    -\nabla \cdot (\phi D \nabla p) - p_N |\nabla\phi| = \phi f.
\end{equation}
The pressure correction from \cref{eq:HIPC2} satisfies
a homogeneous Neumann boundary condition on $\domainbnd{}_D$,
the term $p_N |\nabla\phi|$ disappears from the \DD{} formulation.
However,
boundary conditions that lie on the exterior boundary $\extdomainbnd{}$,
such as the non-extended boundaries,
are enforced directly.
Therefore, as the Neumann boundary condition uses the weak approach,
we apply that approach to the divergence term.
\begin{equation}
    0
    = \int_{\domain{}}{\nabla \cdot u q \; dx}
    = \int_{\domain{}}{ u \cdot \nabla q \; dx} - \int_{\domainbnd{}_D}{ g \cdot n q \; dS}
    \approx \int_{\extdomain{}}{ \phi u \cdot \nabla q + g \cdot \nabla\phi q \;dx}.
\end{equation}
This is implemented in \ddfem{} as an advection term for the pressure,
and using the DDM1 method.

\subsubsection{Taylor-Green}

\label{sec:tg}

The Taylor-Green vortex is a manufactured solution for Navier-Stokes,
so we have an exact solution \citep{Taylor1937}.
Following the example by \citet{Loy2017}, we have the exact solution
\begin{subequations}
    \begin{align}
        u(t, x) & {}= \left(-\cos(2 \pi x_1)\sin(2 \pi x_2) e^{-8 \pi^2 \nu t},\; \sin(2 \pi x_1)\cos(2 \pi x_2) e^{-8 \pi^2 \nu t}\right)^T, \\
        p(t, x) & {}= -\frac{1}{4}(\cos(4\pi x_1)+ \cos(4 \pi x_2))e^{-16 \pi^2 \nu t}.
    \end{align}
\end{subequations}
This gives $f=0$, and we will select $\nu = 0.01$ to solve on the domain
of a ball with a 0.5 radius, $\domain{} = B_{0.5}(0)$,
with Dirichlet boundary conditions using the exact solution.
The problem is solved for $t \in [0, 1]$
using a time step of $\tau = 0.002$.
Numerical experiments show that reductions in the time step do not reduce the error.

For this problem, we will look at the errors across time using the following errors \citep{Sutton2018},
\begin{subequations}
    \begin{align}
        E_{L^2}(u_h)
         & = \left(
            \sum_{n}{ \tau {\|\chi_{\domain{}} (u(t_n)-u^n_h)\|}^2_{L^{2}} }
        \right)^{\frac{1}{2}},
        \\
        E_{H^1}(u_h)
         & =  E_{L^2}(u_h) + \left(
            \sum_{n}{ \tau {\| \chi_{\domain{}} \left( \nabla u(t_n) - \nabla u^n_h \right) \|}^2_{L^2} }
        \right)^{\frac{1}{2}},
        \\
        E_{L^2}(p_h)
         & = \left(
            \sum_{n}{ \tau {\|\chi_{\domain{}} (p(t_n)-p^n_h)\|}^2_{L^{2}} }
        \right)^{\frac{1}{2}}.
    \end{align}
\end{subequations}
Following from \cref{eq:errors},
the mixed methods use $\sigma^n_h$ instead of $\nabla u^n_h$.

The results for \NSDDM{} are shown in \cref{tab:nstg-nsddm} and \MixedZ{} are shown in \cref{tab:nstg-mixed}.
Recall, because this problem has
Dirichlet boundary conditions for the velocity on the whole boundary,
the pressure is solved uniquely up to a constant.
This is addressed by shifting to match the mean of the exact solution,
\begin{equation}
    \label{eq:pressure_shift}
    \frac{1}{|\domain{}|}\int_{\domain{}} p(t_n) - p^n_h \; dx,
\end{equation}
which provides the closest $L^2$ match.
The results show that the magnitude of the errors between the methods is the same,
with \NSDDM{} showing a higher convergence rate.

\begin{table}[htbp]
    \centering
    \begin {tabular}{ccccccc}
    \toprule $h$&$E_{L^2}(u_h)$&$E_{L^2}(u_h)$ EOC&$E_{H^1}(u_h)$&$E_{H^1}(u_h)$ EOC&$E_{L^2}(p_h)$&$E_{L^2}(p_h)$ EOC\\\midrule
    \ensuremath {5.00\cdot 10^{-2}}&\ensuremath {2.07\cdot 10^{-2}}&\ensuremath {-}&\ensuremath {5.45\cdot 10^{-1}}&\ensuremath {-}&\ensuremath {8.59\cdot 10^{-3}}&\ensuremath {-}\\
    \ensuremath {2.50\cdot 10^{-2}}&\ensuremath {8.40\cdot 10^{-3}}&\ensuremath {1.302}&\ensuremath {3.87\cdot 10^{-1}}&\ensuremath {0.495}&\ensuremath {2.93\cdot 10^{-3}}&\ensuremath {1.551}\\
    \ensuremath {1.25\cdot 10^{-2}}&\ensuremath {3.36\cdot 10^{-3}}&\ensuremath {1.322}&\ensuremath {2.84\cdot 10^{-1}}&\ensuremath {0.445}&\ensuremath {8.74\cdot 10^{-4}}&\ensuremath {1.745}\\
    \ensuremath {6.25\cdot 10^{-3}}&\ensuremath {1.39\cdot 10^{-3}}&\ensuremath {1.277}&\ensuremath {2.08\cdot 10^{-1}}&\ensuremath {0.452}&\ensuremath {2.83\cdot 10^{-4}}&\ensuremath {1.630}\\\bottomrule
    \end {tabular}
    \caption{\NSDDM{} method error for Taylor-Green problem.}
    \label{tab:nstg-nsddm}
\end{table}
\begin{table}[htbp]
    \centering
    \begin {tabular}{ccccccc}
    \toprule $h$&$E_{L^2}(u_h)$&$E_{L^2}(u_h)$ EOC&$E_{H^1}(u_h)$&$E_{H^1}(u_h)$ EOC&$E_{L^2}(p_h)$&$E_{L^2}(p_h)$ EOC\\\midrule
    \ensuremath {5.00\cdot 10^{-2}}&\ensuremath {2.09\cdot 10^{-2}}&\ensuremath {-}&\ensuremath {2.47\cdot 10^{-1}}&\ensuremath {-}&\ensuremath {6.69\cdot 10^{-3}}&\ensuremath {-}\\
    \ensuremath {2.50\cdot 10^{-2}}&\ensuremath {8.22\cdot 10^{-3}}&\ensuremath {1.348}&\ensuremath {9.87\cdot 10^{-2}}&\ensuremath {1.323}&\ensuremath {2.31\cdot 10^{-3}}&\ensuremath {1.534}\\
    \ensuremath {1.25\cdot 10^{-2}}&\ensuremath {3.32\cdot 10^{-3}}&\ensuremath {1.307}&\ensuremath {4.25\cdot 10^{-2}}&\ensuremath {1.217}&\ensuremath {7.26\cdot 10^{-4}}&\ensuremath {1.669}\\
    \ensuremath {6.25\cdot 10^{-3}}&\ensuremath {1.42\cdot 10^{-3}}&\ensuremath {1.228}&\ensuremath {2.37\cdot 10^{-2}}&\ensuremath {0.839}&\ensuremath {2.86\cdot 10^{-4}}&\ensuremath {1.346}\\\bottomrule
    \end {tabular}
    \caption{\MixedZ{} method error for Taylor-Green problem.}
    \label{tab:nstg-mixed}
\end{table}

We can see both methods have the velocity converging to the correct solution
and have similar error magnitudes to the previous experiments in \cref{subsection:ad_numerical}.
The $L^2$ EOCs are similar to those in the previous experiments, with \MixedZ{} performing close to \NSDDM{};
although both EOCs fall short of the theoretical rate for Taylor-Green. %
The most significant difference appears in the $H^1$ error,
with \MixedZ{} demonstrating higher EOCs and lower absolute errors.
Furthermore, the pressure error continues to decrease under mesh refinement.
The observed EOCs converge to either $1$ or $0.5$,
consistent with the behaviour shown in \cref{fig:adArc},
suggesting the \DD{} error is dominating.
The trade-off between accuracy and stability,
dependent on the ratio of interfacial width to mesh size ($\epsilon/h$),
and the impact of the Neumann boundary condition approximation,
are left for future work.

\subsubsection{Flow around a cylinder}

\begin{figure}[htbp]
    \centering
    \includegraphics{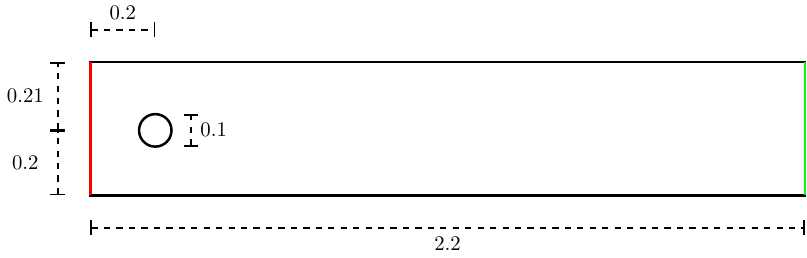}{}
    \caption{Cylinder Flow Domain}
    \label{fig:cyclinder}
\end{figure}

We focus on the Kármán vortex street case from the DFG 2D-3 benchmark \citep{Schaefer1996}.
This problem considers the flow through a domain,
\begin{equation}
    \Omega = [0, 2.2]\times[0, 0.41] \backslash B_{0.05}(0.2, 0.2),
\end{equation}
which is shown in \cref{fig:cyclinder}.
For the \DD{} approach we will extend the domain to the full rectangle,
$B = [0, 2.2]\times[0, 0.41]$.

The problem has $\nu=0.001$, $t\in[0,8]$,
and the following boundary conditions:
the time dependent Dirichlet inflow boundary
\begin{equation}
    u(t, (0, x_2)) = \left(\frac{6 x_2 (0.41 - x_2)}{0.41^2}\sin\left(\frac{\pi t}{8}\right), 0 \right)^T
    \quad \forall x_2 \in [0, 0.41],
\end{equation}
this is the red boundary in \cref{fig:cyclinder},
a do-nothing outflow boundary
\begin{equation}
    (\nu \nabla u - pI) n = 0,
\end{equation}
shown in green,
and the other boundaries have a no slip boundary
\begin{equation}
    u = 0.
\end{equation}
The do-nothing boundary condition is the exterior boundary condition
and far away from the extended region,
therefore, it will not have an impact on the \DD{} approximation.

During the experiments we measure three quantities of interest \citep{Kilian2002,John2002}:
drag coefficient, $C_d(t)$;
lift coefficient, $C_l(t)$;
and the pressure difference around the cylinder,
\begin{equation}
    \Delta p(t) = p(t, (0.15, 0.2)) - p(t, (0.25,0.2)).
\end{equation}
These coefficients are derived from
\begin{equation}
    \begin{pmatrix}
        C_d \\ C_l
    \end{pmatrix}
    =
    \frac{2}{\rho L U^2_{\text{mean}}}
    \int_{\domainbnd{}}{ (\nu \nabla u - p I) n \;dS}.
\end{equation}
This can be easily converted into an integral that we can compute in the \DD{} framework;
we have the extended normal defined in \cref{eq:normal} and the surface delta in \cref{eq:surfacedelta}.
\begin{equation}
    \begin{pmatrix}
        C_d \\ C_l
    \end{pmatrix}
    \approx
    \frac{2}{\rho L U^2_{\text{mean}}}
    \int_{\extdomain{}}{ (p I - \nu \nabla u) \nabla \phi \;dx}.
\end{equation}

We use a refined mesh such that
the largest elements have a diameter of $h_{\text{max}}=0.02$
and the diameter of the smallest elements is $h_{\text{min}}=0.005$.
The element size increases linearly,
starting at a distance of $5\epsilon = 0.0875$
from the interface of the cylinder,
reaching the largest size at $20\epsilon=0.35$.

\begin{figure}[htbp]
    \centering
    \includegraphics{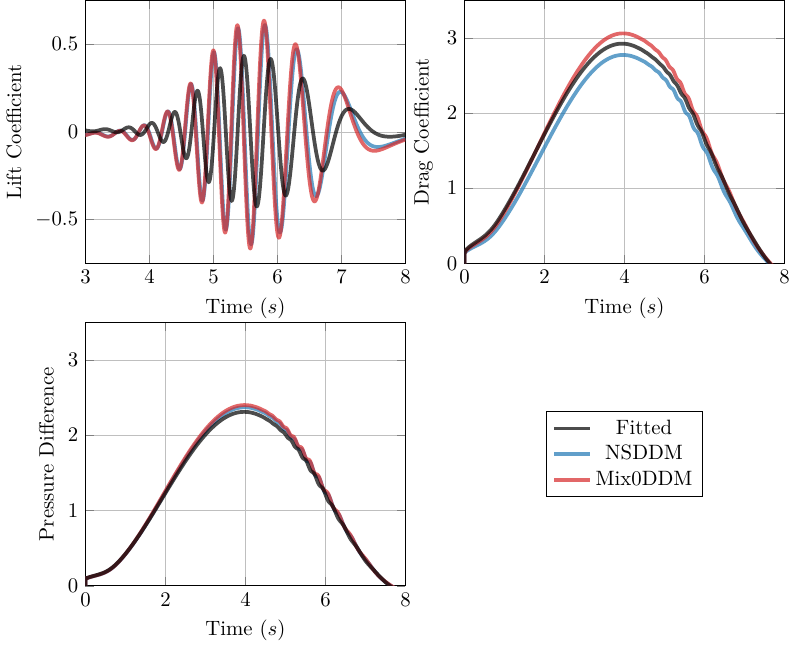}{}
    \caption{Measurements for Navier-Stokes flow around cylinder.}
    \label{fig:nscoef}
\end{figure}

To ensure numerical stability and reduce time step errors,
we set $\tau = 0.005$.
Results are shown in \cref{fig:nscoef},
alongside a simulation using the same IPC method on a similarly refined, fitted mesh.
These results demonstrate
that our \DD{} framework successfully captures the unsteady physics of the Kármán vortex shedding.
The timing of the drag coefficient and lift coefficient oscillations closely matches the fitted solution.
However, the measured values for these coefficients show a deviation from the fitted solution.

To evaluate the accuracy of these results,
they are compared in \cref{tab:nscoef} against the experiments performed by \citet{John2004}.
\begin{table}[htbp]
    \centering
    \begin {tabular}{lccccc}
    \toprule Method&$t(C_{d,\max })$&$C_{d,\max }$&$t(C_{l,\max })$&$C_{l,\max }$&$\Delta p(8\text {s})$\\\midrule
    \citet {John2004}&\ensuremath {3.93625}&\ensuremath {2.950921575}&\ensuremath {5.693125}&\ensuremath {0.4779500}&\ensuremath {-0.1116000}\\
    Fitted&\ensuremath {3.93500}&\ensuremath {2.924793767}&\ensuremath {5.475000}&\ensuremath {0.4328733}&\ensuremath {-0.1107880}\\
    \MixedZ {}&\ensuremath {3.95600}&\ensuremath {3.061619554}&\ensuremath {5.576000}&\ensuremath {0.6633196}&\ensuremath {-0.0947517}\\
    \NSDDM {}&\ensuremath {3.95600}&\ensuremath {2.774583223}&\ensuremath {5.588000}&\ensuremath {0.6392941}&\ensuremath {-0.0935027}\\\bottomrule
    \end {tabular}
    \caption{Comparison with reference value of flow around a cylinder.}
    \label{tab:nscoef}
\end{table}
From the lift coefficient measurements,
it is clear both methods produce the expected vortex shedding.
Furthermore, the maxima occur close to the reference time; however, both methods overestimate the maximum lift,
with \NSDDM{} producing a value closer to the reference.

A similar behaviour can be spotted in the drag coefficient,
with both methods producing the maximum value close to the expected time.
However, both methods have a noticeable difference from the expected drag coefficient.
The maximum drag for \NSDDM{} undershot the reference value,
and overshot the reference value for \MixedZ{}.
Finally, comparing the pressure difference we can see
both methods produce very close results on the figure,
which is similar to the performance shown for the Taylor-Green problem.

This discrepancy highlights the fundamental compromise of the \DD{} method.
Approximating the sharp boundary with a smooth, diffuse interface,
causes a compromise in accuracy for the precise values.
For many applications the complex geometry is the primary challenge
and this level of accuracy is sufficient.
The decision between a \DD{} approach and a traditional fitted method
depends on the specific goals of the simulation.

Furthermore,
experiments with the same minimum mesh size
but varying $\epsilon$
demonstrate the impact of the diffuse interface width
on the accuracy of the results,
as shown in \cref{tab:ns-eps}.
Note that the refinement region around the cylinder depends on $\epsilon$,
so a larger $\epsilon$ results in more DoFs.
When the ratio of $\epsilon$ to $h$ is too small,
the simulations fail.
This is likely because
the diffuse interface is not sufficiently resolved by the mesh.
For larger stable ratios,
increasing $\epsilon$ reduces the accuracy of the computed drag and lift coefficients
and the pressure difference compared to the reference values.
This occurs because the interface becomes too wide
and impacts the bulk solution. In our simulations we have chosen $W=7$
which our tests indicate provides a good compromise of accuracy and stability.
\begin{table}[htbp]
    \centering
    \begin {tabular}{ccccccc}
    \toprule $\epsilon $&$W$&$t(C_{d,\max })$&$C_{d,\max }$&$t(C_{l,\max })$&$C_{l,\max }$&$\Delta p(8\text {s})$\\\midrule
    \ensuremath {7.50\cdot 10^{-3}}&\ensuremath {3}&\ensuremath {\hbox {NaN}}&\ensuremath {\hbox {NaN}}&\ensuremath {\hbox {NaN}}&\ensuremath {\hbox {NaN}}&\ensuremath {\hbox {NaN}}\\
    \ensuremath {1.25\cdot 10^{-2}}&\ensuremath {5}&\ensuremath {3.95000}&\ensuremath {2.989118059}&\ensuremath {5.548000}&\ensuremath {0.6150220}&\ensuremath {-0.1004733}\\
    \ensuremath {1.75\cdot 10^{-2}}&\ensuremath {7}&\ensuremath {3.95600}&\ensuremath {3.061619554}&\ensuremath {5.576000}&\ensuremath {0.6633196}&\ensuremath {-0.0947517}\\
    \ensuremath {2.25\cdot 10^{-2}}&\ensuremath {9}&\ensuremath {3.96200}&\ensuremath {3.139127942}&\ensuremath {5.604000}&\ensuremath {0.6966963}&\ensuremath {-0.0938075}\\
    \ensuremath {2.75\cdot 10^{-2}}&\ensuremath {11}&\ensuremath {3.96800}&\ensuremath {3.234007875}&\ensuremath {5.636000}&\ensuremath {0.7212485}&\ensuremath {-0.0971956}\\\bottomrule
    \end {tabular}
    \caption{Measurements of flow around a cylinder with \MixedZ{} for varied $\epsilon$.
        Here, $W=\frac{2\epsilon}{h_{\text{min}}}$ is the number of elements across the interfacial region.}
    \label{tab:ns-eps}
\end{table}

Finally, \cref{tab:ns-h} shows the results for \MixedZ{} under mesh
refinement.
As expected, the results improve when is the element diameter is reduced,
with the drag and lift coefficients and the pressure difference
approaching the reference values. Only the timing for the maximum lift
is moving towards a different value which indicates that some further
investigations should be carried out which are beyond the scope of this
paper.
\begin{table}[htbp]
    \centering
    \begin {tabular}{llcccccc}
    \toprule $h_{\max }$&$h_{\min }$&$\epsilon $&$t(C_{d,\max })$&$C_{d,\max }$&$t(C_{l,\max })$&$C_{l,\max }$&$\Delta p(8\text {s})$\\\midrule
    0.01&0.0025&\ensuremath {8.75\cdot 10^{-3}}&\ensuremath {3.94400}&\ensuremath {2.940791521}&\ensuremath {5.526000}&\ensuremath {0.5697826}&\ensuremath {-0.1067541}\\
    0.02&0.005&\ensuremath {1.75\cdot 10^{-2}}&\ensuremath {3.95600}&\ensuremath {3.061619554}&\ensuremath {5.576000}&\ensuremath {0.6633196}&\ensuremath {-0.0947517}\\
    0.04&0.01&\ensuremath {3.50\cdot 10^{-2}}&\ensuremath {3.97600}&\ensuremath {3.466204277}&\ensuremath {5.690000}&\ensuremath {0.7707407}&\ensuremath {-0.0982857}\\\bottomrule
    \multicolumn {3}{c}{\citet {John2004}} & 3.936 & 2.950921575 & 5.693 & 0.47795 & -0.1116\\ \end {tabular}
    \caption{Measurements of flow around a cylinder with \MixedZ{} for varied mesh size.}
    \label{tab:ns-h}
\end{table}

\section{Conclusion}
\label{section:conclusion}

We have introduced and analysed a new approach for
applying the \DD{} method to partial differential equations on complex domains,
derived from a mixed formulation.
This led to the development of two new methods,
\MixedZ{} and \MixedO{}.
A key benefit of these methods, for certain problems, is that
they transform the essential Dirichlet boundary conditions
into natural boundary conditions
within the gradient equation of the mixed system.
Furthermore, we derived new stabilised methods derived from Nitsche's method,
\NSDDM{} and \NDDM{},
and improved stabilisation of existing methods with a focus on advection
dominated problems.

Numerical experiments for Advection-Diffusion problems demonstrated
that the new methods improve on the existing methods.
\NSDDM{} consistently produced the best $L^2$ error,
and while \MixedZ{} performed close in the $L^2$ error,
it performed significantly better in the $H^1$ error.
This suggests that \MixedZ{} is especially beneficial for problems where
the accuracy of the solution's gradient is important,
due to the boundary conditions being included in $\sigma \approx \nabla u$.

Notably, the \MixedZ{} method was the only method proven to be
coercive without an additional advection stabilisation term,
whereas \NSDDM{} is a high performing symmetric formulation.

The experiments of the incompressible Navier-Stokes equations,
using an incremental pressure correction scheme,
demonstrated their accuracy for complex Dirichlet boundary problems.
Both methods produced the expected behaviour,
however, the measurements showed some error relative to the reference values.

Future work should focus on
a rigorous asymptotic analysis of the new methods
to establish their theoretical convergence rates,
and calculate the optimal choice of $\epsilon$ scaling in $\BC{}$ terms.

\bibliographystyle{elsarticle-num-names}
\bibliography{references.bib}

\FloatBarrier

\appendix

\section{Full Advection-Diffusion Numerical Results}
\label{sec:fullnumerical}

\begin{figure}[!htbp]
    \centering
    \includegraphics{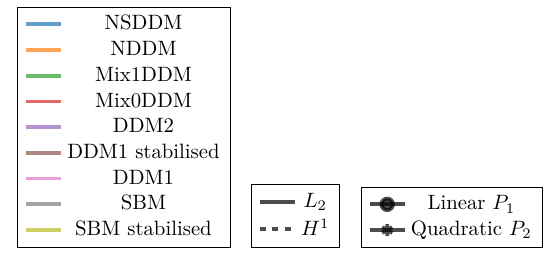}{}
    \caption{Legends for all plots in \ref{sec:fullnumerical}.}
    \label{fig:fulllegends}
\end{figure}

Here we show the results for all approaches tested for
the problems in \cref{subsection:ad_numerical}.
\Cref{fig:fulllegends} shows the legend for the approaches, the
order of basis functions, and the error measurement;
this matches the earlier plots.

\subsection{Diffusion Dominated}

\begin{figure}[!htbp]
    \centering
    \includegraphics{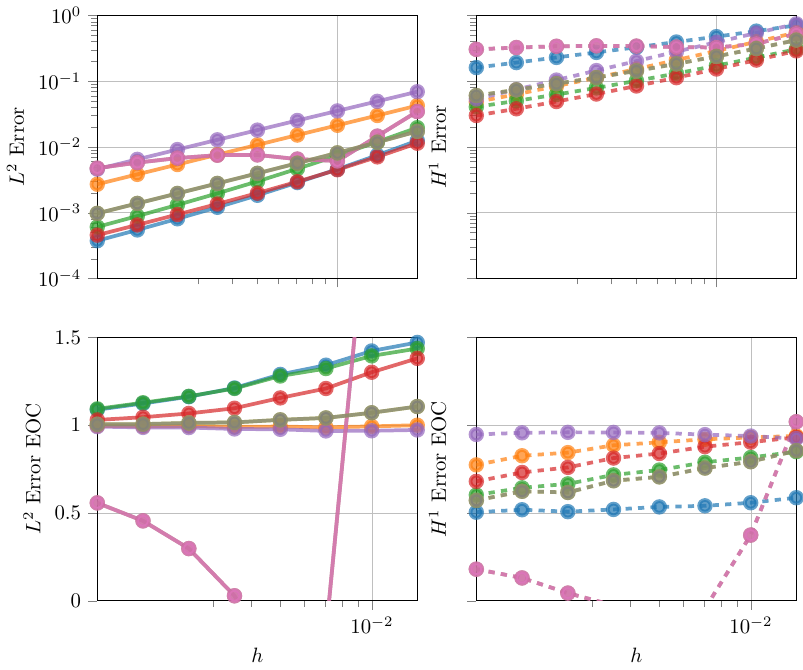}{}
    \caption{Diffusion Dominated problem error plot for Arc domain with linear elements.}
    \label{fig:appDDarc1}
\end{figure}

\begin{figure}[!htbp]
    \centering
    \includegraphics{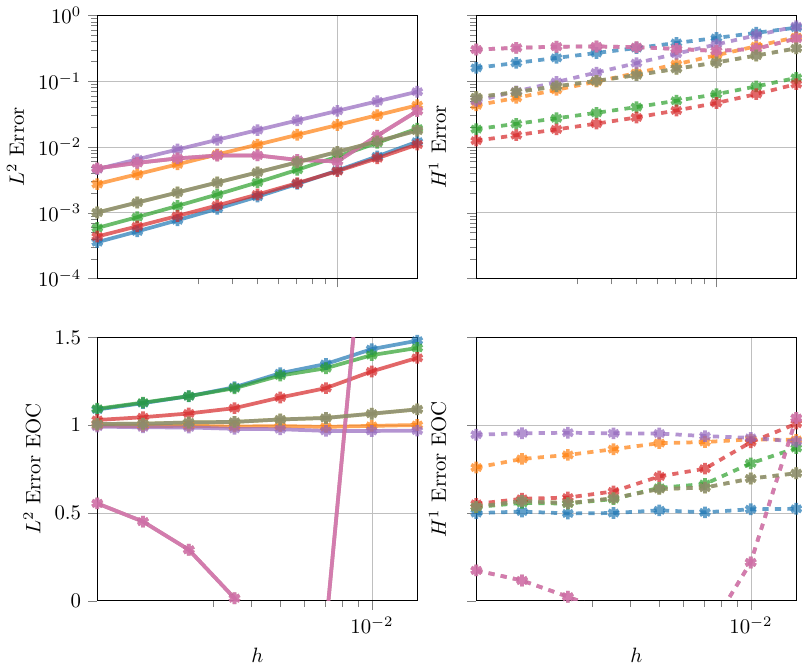}{}
    \caption{Diffusion Dominated problem error plot for Arc domain with quadratic elements.}
    \label{fig:appDDarc2}
\end{figure}

\begin{figure}[!htbp]
    \centering
    \includegraphics{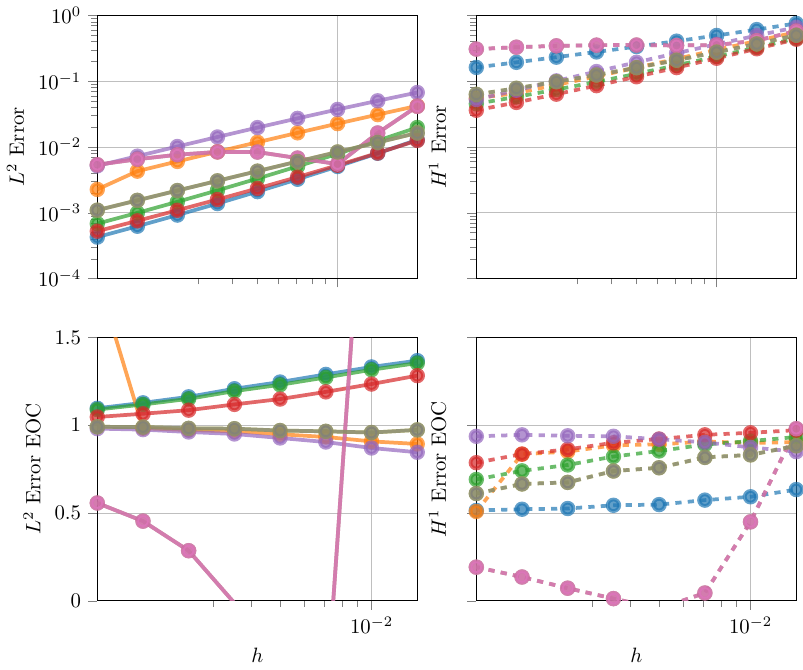}{}
    \caption{Diffusion Dominated problem error plot for Inverted Arc domain with linear elements.}
    \label{fig:appDDiarc1}
\end{figure}

\begin{figure}[!htbp]
    \centering
    \includegraphics{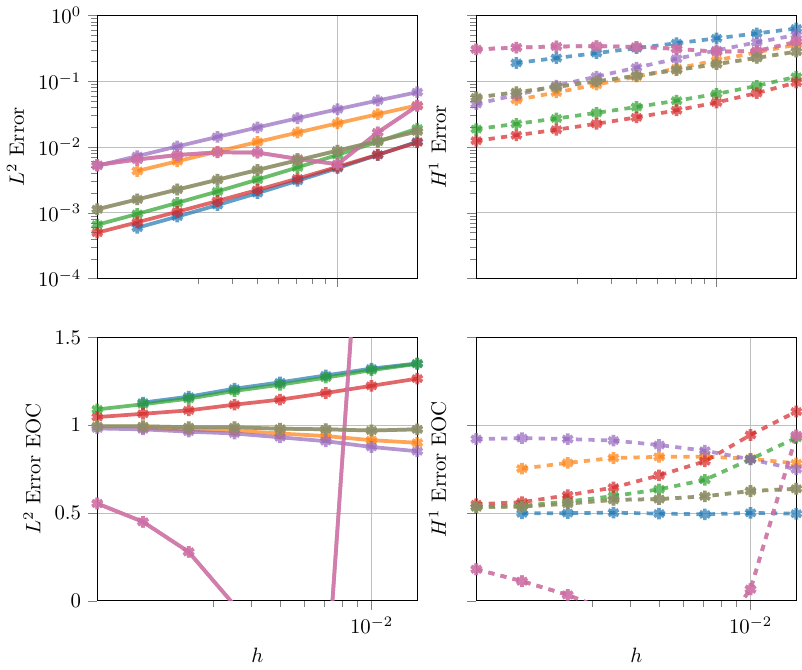}{}
    \caption{Diffusion Dominated problem error plot for Inverted Arc domain with quadratic elements.}
    \label{fig:appDDiarc2}
\end{figure}

\Cref{fig:appDDarc1,fig:appDDarc2,fig:appDDiarc1,fig:appDDiarc2} show
the results for the diffusion dominated problem.
It is clear that \MixedZ{} produces lower errors than \MixedO{} for both the $L^2$ and $H^1$ norms,
and \NSDDM{} achieves a lower $L^2$ error but a larger $H^1$ error than \NDDM{}.
Importantly, they are showing lower errors and better convergence than
the traditional approaches of DDM1 and SBM.

\FloatBarrier
\subsection{Advection Dominated}

\begin{figure}[!htbp]
    \centering
    \includegraphics{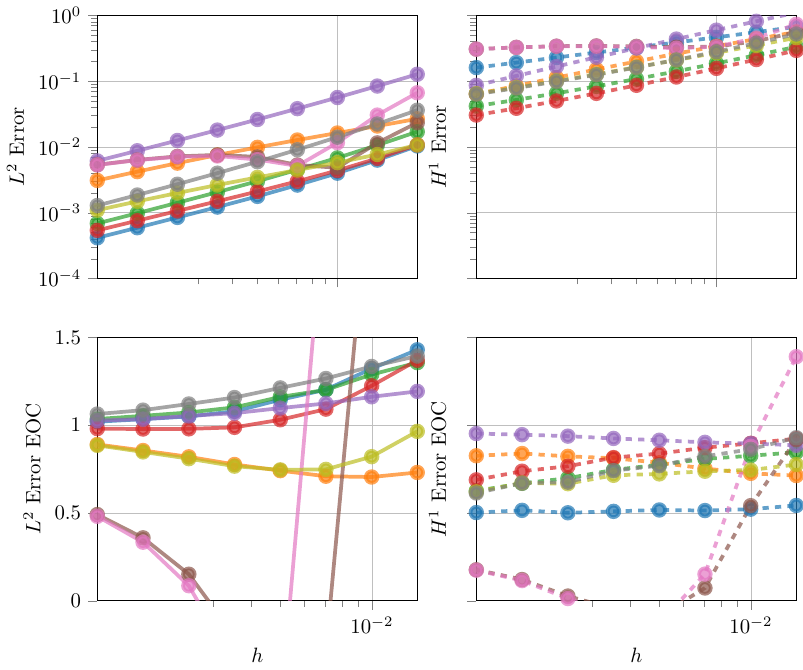}{}
    \caption{Advection Dominated problem error plot for Arc domain with linear elements.}
    \label{fig:appADarc1}
\end{figure}

\begin{figure}[!htbp]
    \centering
    \includegraphics{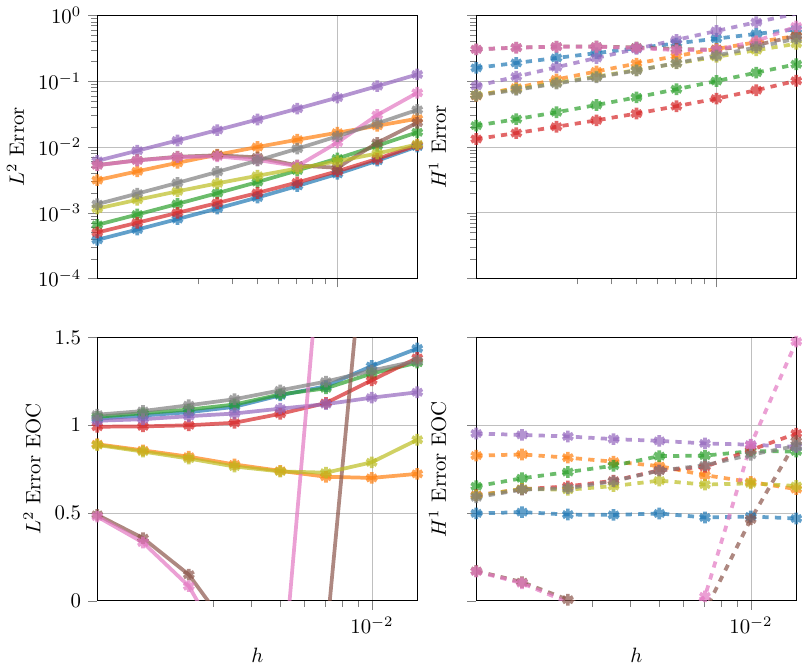}{}
    \caption{Advection Dominated problem error plot for Arc domain with quadratic elements.}
    \label{fig:appADarc2}
\end{figure}

\begin{figure}[!htbp]
    \centering
    \includegraphics{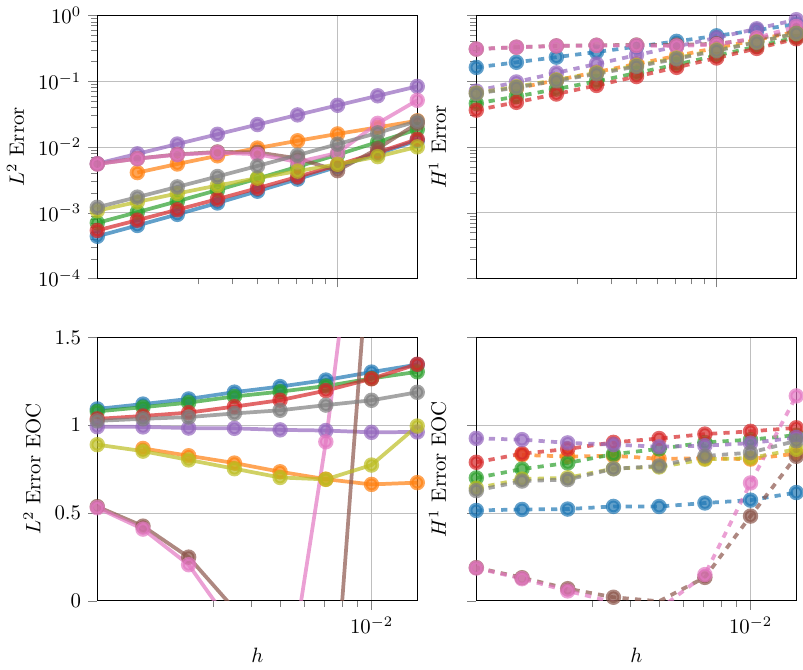}{}
    \caption{Advection Dominated problem error plot for Inverted Arc domain with linear elements.}
    \label{fig:appADiarc1}
\end{figure}

\begin{figure}[!htbp]
    \centering
    \includegraphics{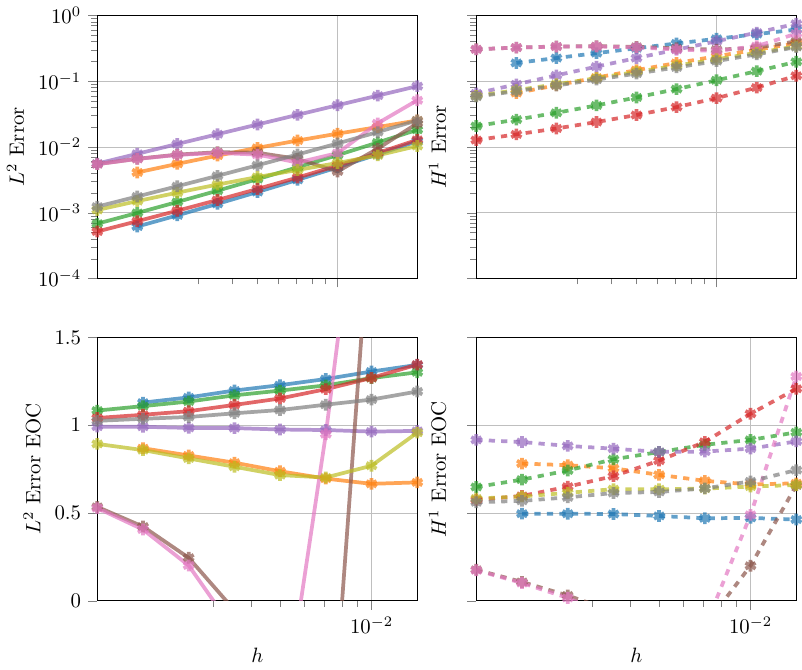}{}
    \caption{Advection Dominated problem error plot for Inverted Arc domain with quadratic elements.}
    \label{fig:appADiarc2}
\end{figure}

\Cref{fig:appADarc1,fig:appADarc2,fig:appADiarc1,fig:appADiarc2} show
the results for the advection dominated problem,
which are similar to the above diffusion dominated problem.
We can see the impact of the stabilised terms on DDM1 and SBM,
which show the stabilised SBM error has a smaller magnitude but slower convergence.

\FloatBarrier

\section{Reconstructed Gradient Error Comparison}
\label{sec:gradError}

\begin{figure}[htbp]
    \centering
    \begin{subfigure}[t]{\textwidth}
        \centering
        \includegraphics[width=\linewidth]{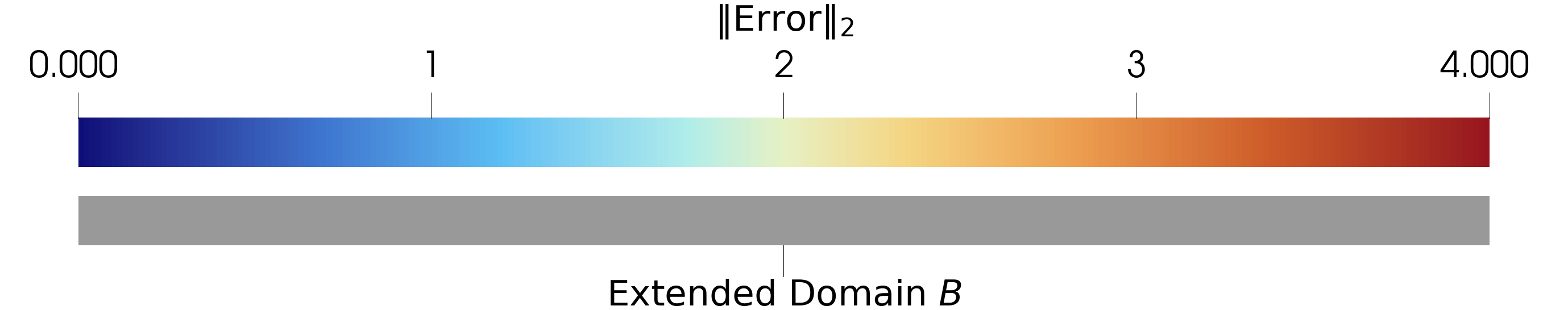}
    \end{subfigure}
    \hfill
    \begin{subfigure}{0.495\textwidth}
        \centering
        \includegraphics[width=\linewidth]{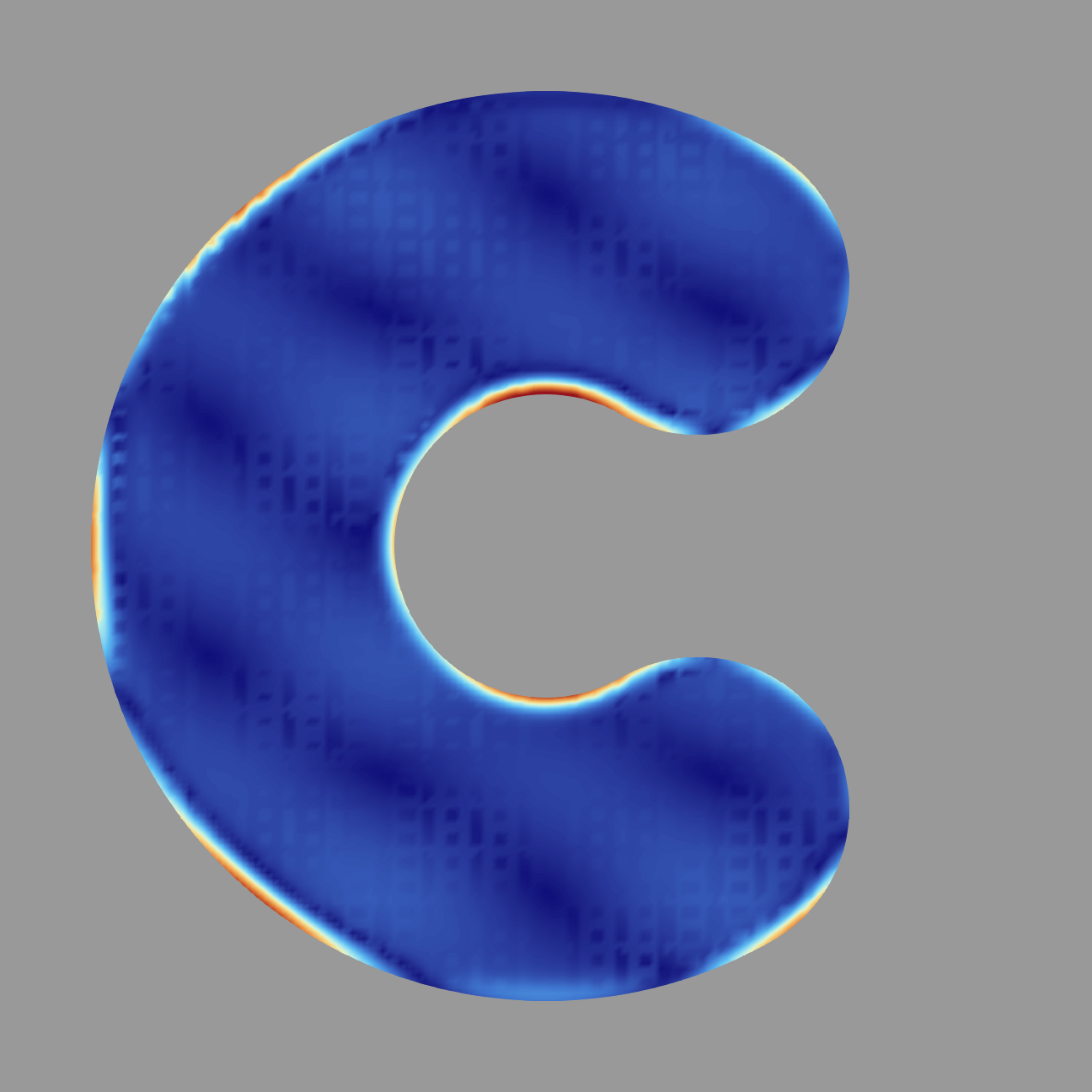}
        \caption{$|\nabla u - \nabla u_h|$.}
    \end{subfigure}
    \hfill
    \begin{subfigure}{0.495\textwidth}
        \centering
        \includegraphics[width=\linewidth]{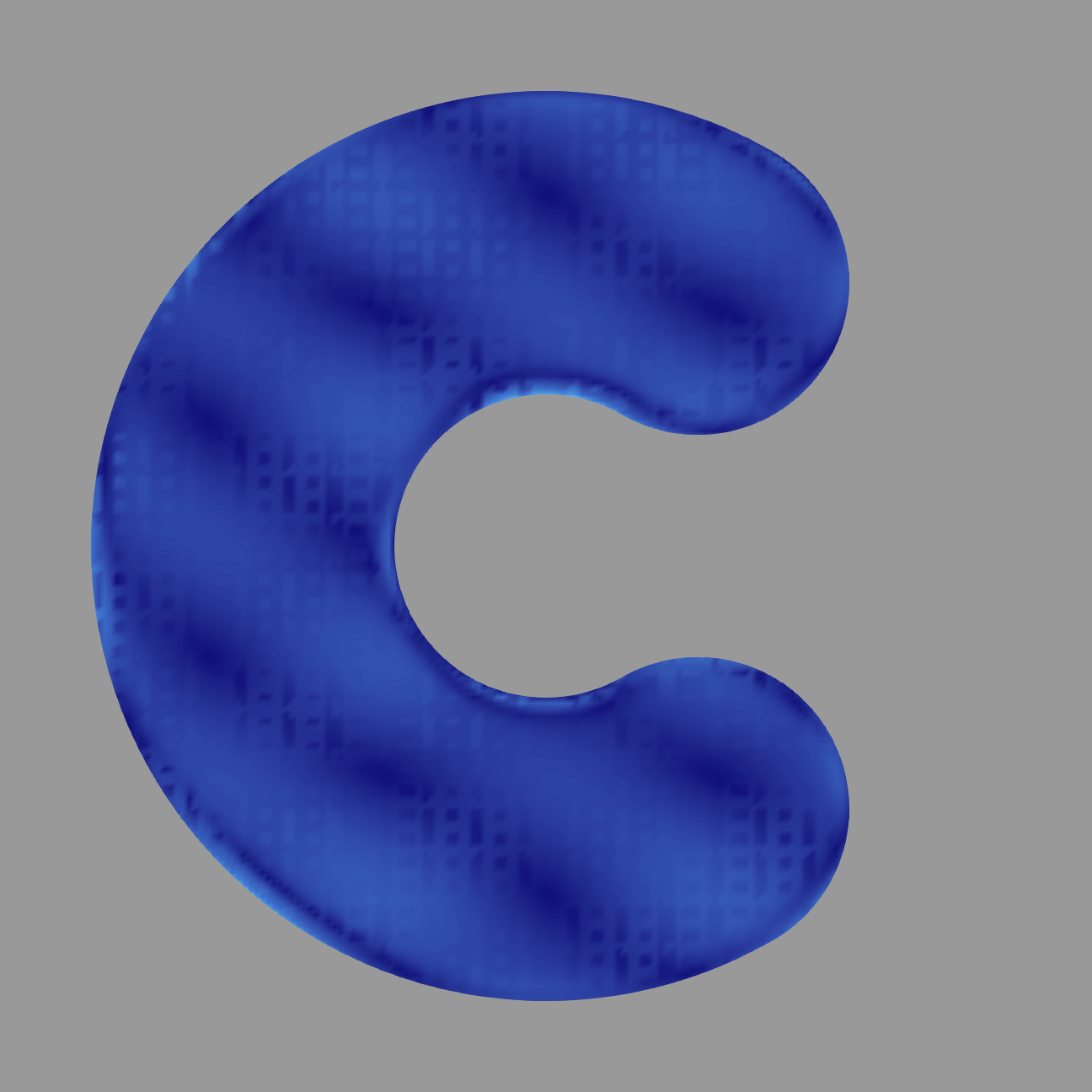}
        \caption{$|\nabla u - \sigma_h|$.}
    \end{subfigure}
    \hfill
    \caption{Gradient errors for \MixedZ{} with linear elements for the diffusion problem on the arc domain.}
    \label{fig:gradientErrorComparison}
\end{figure}

To show the effect of the gradient reconstruction used in \cref{eq:errors},
we consider the diffusion problem in \cref{subsection:ad_numerical} on the arc domain.
\Cref{fig:gradientErrorComparison} compares pointwise the gradient errors for \MixedZ{}
with linear elements at the second mesh refinement.
The reconstructed gradient $\sigma_h$ reduces the error near the boundary
compared to direct gradient $\nabla u_h$,
while the errors in the interior are similar.
This motivates the improved boundary accuracy discussed in \cref{sec:Numerical}
for the mixed methods.

\end{document}